\documentclass[11pt]{article}%
\usepackage{amsmath}
\usepackage{graphicx}
\usepackage{cite}
\usepackage{subfigure}
\usepackage{multirow}
\usepackage{authblk}
\usepackage{indentfirst}
\usepackage{caption}
\usepackage{algorithm}
\usepackage{algorithmicx}
\usepackage{algpseudocode}
\usepackage{graphicx}
\usepackage{bm}
\usepackage{CJK}
\usepackage{geometry}
\usepackage{float}
\usepackage{colortbl}
\usepackage{textcomp,booktabs}
\usepackage{array}
\usepackage{xcolor}%
\usepackage{amsfonts}%
\usepackage{amssymb}

\newcommand{\eqnb}{\begin{equation}}
\newcommand{\eqne}{\end{equation}}

\begin{document}

\title{An Overview for Markov Decision Processes in Queues and Networks}
\author{Quan-Lin Li$^{a}$, Jing-Yu Ma$^{b}$, Rui-Na Fan$^{b}$, Li Xia$^{c}$\\$^{a}$School of Economics and Management,\\Beijing University of Technology, Beijing 100124, China\\$^{b}$School of Economics and Management,\\Yanshan University, Qinhuangdao 066004, China\\$^{c}$Bussiness School, Sun Yat-sen University, Guangzhou 510275, China}
\maketitle

\begin{abstract}
Markov decision processes (MDPs) in queues and networks have been an
interesting topic in many practical areas since the 1960s. This paper Provides
a detailed overview on this topic and tracks the evolution of many basic
results. Also, this paper summarizes several interesting directions in the
future research. We hope that this overview can shed light to MDPs in queues
and networks, and also to their extensive applications in various practical areas.

\vskip                 0.5cm


\end{abstract}

\section{Introduction}

One main purpose of this paper is to provide an overview for research on MDPs
in queues and networks in the last six decades. Also, such a survey is first
related to several other basic studies, such as, Markov processes, queueing
systems, queueing networks, Markov decision processes, sensitivity-based
optimization, stochastic optimization, fluid and diffusion control. Therefore,
our analysis begins from three simple introductions: Markov processes and
Markov decision processes, queues and queueing networks, and queueing dynamic control.

\textbf{(a) Markov processes and Markov decision processes}

The Markov processes, together with the Markov property, were first introduced
by a Russian mathematician: Andrei Andreevich Markov (1856-1922) in 1906. See
Markov \cite{Mar:1906} for more details. From then on, as a basically
mathematical tool, the Markov processes have extensively been discussed by
many authors, e.g., see some excellent books by Doob \cite{Doo:1953}, Karlin
\cite{Kar:1968}, Karlin and Taylor \cite{Kar:1981}, Chung \cite{Chu:1967},
Anderson \cite{And:1991}, Kemeny \textit{et al}. \cite{Kem:1976}, Meyn and
Tweedie \cite{Mey:1996}, Chen \cite{Chen:2004}, Ethier and Kurtz
\cite{Eth:2005} and so on.

In 1960, Howard \cite{How:1960} is the first to propose and discuss the MDP
(or stochastic dynamic programming) in terms of his Ph.D thesis, which opened
up a new and important field through an interesting intersection between
Markov processes and dynamic programming (e.g., see Bellman and Kalaba
\cite{Bel:1965}). From then on, not only are the MDPs an important branch in
the area of Markov processes, but also it is a basic method in modern dynamic
control theory. Crucially, the MDPs have been greatly motivated and widely
applied in many practical areas in the past 60 years. Readers may refer to
some excellent books, for example, the discrete-time MDPs by Puterman
\cite{Pur:1994}, Glasserman and Yao \cite{Gla:1994}, Bertsekas \cite{Ber:1996}%
, Bertsekas and Tsitsiklis \cite{BerT:1996}, Hern\'{a}dez-Lerma and Lasserre
\cite{Her:1996, Her:1999}, Altman \cite{Alt:1999}, Koole \cite{Koo:2007} and
Hu and Yue \cite{Hu:2007}; the continuous-time MDPs by Guo and
Hern\'{a}ndez-Lerma \cite{Guo:2009}; the partially observable MDPs by
Cassandra \cite{Cas:1998} and Krishnamurthy \cite{Kri:2016}; the competitive
MDPs (i.e., stochastic game) by \cite{Fil:2012}; the sensitivity-based
optimization by Cao \cite{Cao:2007}; some applications of MDPs by Feinberg and
Shwartz (Eds.) \cite{FeiS:2002}\ and Boucherie and Van Dijk (Eds.)
\cite{Bou:2017}; and so on.

\textbf{(b) Queues and queueing networks}

In the early 20th century, a Danmark mathematician: Agner Krarup Erlang,
published a pioneering work \cite{Erl:1909} of queueing theory in 1909, which
started the study of queueing theory and traffic engineering. Over the past
100 years, queueing theory has been regarded as a key mathematical tool not
only for analyzing practical stochastic systems but also for promoting theory
of stochastic processes (such as Markov processes, semi-Markov processes,
Markov renew processes, random walks, martingale theory, fluid and diffusion
approximation, and stochastic differential equations). On the other hand, the
theory of stochastic processes can support and carry forward advances in
queueing theory and applications (for example single-server queues,
multi-server queues, tandem queues, parallel queues, fork-join queues, and
queueing networks). It is worthwhile to note that so far queueing theory has
been widely applied in many practical areas, such as manufacturing systems,
computer and communication networks, transportation networks, service
management, supply chain management, sharing economics, healthcare and so forth.

\textit{The single-server queues and the multi-server queues:} In the early
development of queueing theory (1910s to 1970s), the single-server queues were
a main topic with key results including Khintchine formula, Little's law,
birth-death processes of Markovian queues, the embedded Markov chain, the
supplementary variable method, the complex function method and so on. In 1969,
Professor J.W. Cohen published a wonderful summative book \cite{Coh:1969} with
respect to theoretical progress of single-server queues.

It is a key advance that Professor M.F. Neuts proposed and developed the
phase-type (PH) distributions, Markovian arrival processes (MAPs), and the
matrix-geometric solution, which were developed as the matrix-analytic method
in the later study, e.g., see Neuts \cite{Neu:1981, Neu:1989} and and Latouche
and Ramaswami \cite{Lat:1999} for more details. Further, Li \cite{Li:2010}
proposed and developed the RG-factorizations for any generally irreducible
block-structured Markov processes. Crucially, the RG-factorizations promote
the matrix-analytic method to a unified matrix framework both for the
steady-state solution and for the transient solution (for instance the first
passage time and the sojourn time). In addition, the matrix-analytic method
and the RG-factorizations can effectively deal with small-scale stochastic
models with several nodes.

In the study of queueing systems, some excellent books include Kleinrock
\cite{Kle:1975, Kle:1976}, Tijms \cite{Tij:1994} and Asmussen \cite{Asm:2003}.
Also, an excellent survey on key queueing advances was given in Syski
\cite{Sys:1997}; and some overview papers on different research directions
were reported by top queueing experts in two interesting books by Dshalalow
\cite{Dsh:1995, Dsh:1997}.

\textit{The queueing networks: }In 1957, J.R. Jackson published a seminal
paper \cite{Jac:1957} which started research on queueing networks. Subsequent
interesting results include Jackson \cite{Jac:1963}, Baskett et al.
\cite{Bas:1975}, Kelly \cite{Kel:1976, Kel:1991}, Disney and K\"{o}nig
\cite{Dis:1985}, Dobrushin \textit{et al}. \cite{Dob:1990}, Harrison
\cite{Har:1985}, Dai \cite{Dai:1995} and so on. For the queueing networks, the
well-known examples contain Jackson networks, BCMP networks, parallel
networks, tandem networks, open networks, closed networks, polling queues,
fork-join networks and distributed networks. Also, the product-form solution,
the quasi-reversibility and some approximation algorithms are the basic
results in the study of queueing networks.

For the queueing networks, we refer readers to some excellent books such as
Kelly \cite{Kel:1979}, Van Dijk \cite{Dij:1993}, Gelenbe \textit{et al}.
\cite{Gel:1998}, Chao \textit{et al}. \cite{Cha:1999}, Serfozo \cite{Ser:1999}%
, Chen and Yao \cite{ChenY:2001}, Balsamo \textit{et al}. \cite{Bal:2001},
Daduna \cite{Dad:2001}, Bolch \textit{et al}. \cite{Bol:2006} and Boucherie
and Van Dijk (Eds.) \cite{Bou:2011}.

For applications of queueing networks, readers may refer to some excellent
books, for example, manufacturing systems by Buzacott and Shanthikumar
\cite{Buz:1993}, communication networks by Chang \cite{Cha:2000}, traffic
networks by Garavello and Piccoli \cite{Gar:2006}, healthcare by Lakshmi and
Iyer \cite{Lak:2013}, service management by Demirkan et al. \cite{Dem:2011}
and others.

\textbf{(c) Queueing dynamic control}

In 1967, Miller \cite{Mil:1967} and Ryokov and Lembert \cite{Ryk:1967} are the
first to apply the MDPs to consider dynamic control of queues and networks.
Those two works opened a novel interesting research direction: MDPs in queues
and networks.

For MDPs of queues and networks, we refer readers to three excellent books by
Kitaev and Rykov \cite{Kit:1995}, Sennott \cite{Sen:2009} and Stidham
\cite{Sti:2009}.

In MDPs of queues and networks, so far there have been some best survey
papers, for instance, Crabill \textit{et al}. \cite{Cra:1973, Cra:1977}, Sobel
\cite{Sob:1974}, Stidham and Prabhu \cite{Sti:1974}, Rykov \cite{Ryk:1975,
Ryk:2017}, Kumar \cite{Kum:1990}, Stidham and Weber \cite{Sti:1993}, Stidham
\cite{Sti:2002} and Brouns \cite{Bro:2003}.

For some Ph.D thesises by using MDPs of queues and networks, reader may refer
to, such as, Farrell \cite{Far:1976}, Abdel-Gawad \cite{Abd:1984}, Bartroli
\cite{Bar:1989}, Farrar \cite{Far:1992}, Veatch \cite{Vea:1992}, Altman
\cite{Alt:1994}, Atan \cite{Ata:1997} and Efrosinin \cite{Efr:2004}.

Now, MDPs of queues and networks play an important role in dynamic control of
many practical stochastic networks, for example, inventory control
\cite{Fed:1984, Fed:1992, Cao:2016}, supply chain management \cite{Ett:2000},
maintenance and quality \cite{Kuo:2006, Dim:2008}, manufacturing systems
\cite{Jo:1991, Buz:1993}, production lines \cite{Wu:2010}, communication
networks \cite{Alt:2002, Paj:2014, Als:2015}, wireless and mobile networks
\cite{Ahm:2005, Din:2013}, cloud service \cite{Sun:2014}, healthcare
\cite{Pat:2011}, airport management \cite{Rue:1985, Lau:1999},
energy-efficient management \cite{Qiu:1999, Oka:2015} and artificial
intelligence \cite{Kol:2012, Sig:2013}. With rapid development of Internet of
Things (IoT), big data, cloud computing, blockchain and artificial
intelligence, it is necessary to discuss MDPs of queues and networks under an
intelligent environment.

\vskip       0.3cm

From the detailed survey on MDPs of queues and networks, this paper suggests a
future research under an intelligent environment from three different levels
as follows:

\begin{enumerate}
\item \textit{Networks with several nodes:} Analyzing MDPs of policy-based
Markov processes with block structure, for example, QBD processes, Markov
processes of GI/M/1 type, and Markov processes of M/G/1 type, and
specifically, discussing their sensitivity-based optimization.

\item \textit{Networks with a lot of nodes:} discussing MDPs of practical big
networks, such as blockchain systems, sharing economics, intelligence
healthcare and so forth.

\item \textit{Networks with a lot of clusters:} studying MDPs of practical big
networks by means of the mean-field theory, e.g., see Gast and Gaujal
\cite{Gas:2011}, Gast et al. \cite{Gas:2012} and Li \cite{Li:2016}.
\end{enumerate}

\vskip       0.3cm

The remainder of this paper is organized as follows. Sections 2 to 5 provide
an overview for MDPs of single-server queues, multi-server queues, queueing
networks, and queueing networks with special structures, respectively. Section
6 sets up specific objectives to provide an overview for key objectives in queueing dynamic control. Section 7 introduce the sensitivity-based optimization and the
event-based optimization, both of which are applied to analyze MDPs of queues
and networks. Finally, we give some concluding remarks in Section 8.

\section{\textbf{MDPs of Single-Server Queues}}

In this section, we provide an overview for MDPs of single-server queues,
including the M/M/1 queues, the M/M/c queues, the M/G/1 queues, the GI/M/1
queues and others. In the early research on MDPs of queues and networks, the
single-server queues have been an active topic for many years.

\textbf{(1) MDPs of M/M/1 queues}

Kofman and Lippman \cite{Kof:1981}, Rue and Rosenshine \cite{Rue:1981a,
Rue:1981b}, Yeh and Thomas \cite{Yeh:1983}, Lu and Serfozo \cite{Lu:1984},
Plum \cite{Plu:1991}, Altman \cite{Alt:1996}, Kitaev and Serfozo
\cite{Kit:1999}, Sava\c{s}aneril et al. \cite{Sav:2010} and Dimitrakopoulos
and Burnetas \cite{Dim:2017}.

\textbf{(2) MDPs of M/G/1 queues}

Mitchell \cite{Mit:1973}, Doshi \cite{Dos:1977, Dos:1978}, Gallisch
\cite{Gal:1979}, Rue and Rosenshine \cite{Rue:1983}, Jo and Stidham
\cite{Jo:1983}, Mandelbaum and Yechiali \cite{Man:1983}, Kella \cite{Kel:1990}%
, Wakuta \cite{Wak:1991}, Altman and Nain \cite{Alt:1993}, Feinberg and Kim
\cite{Fei:1996}, Feinberg and Kella \cite{Fei:2002} and Sanajian et al.
\cite{San:2010}.

\textbf{(3) MDPs of GI/M/1 queues}

Stidham \cite{Sti:1978} and Mendelson and Yechiali \cite{Men:1981}.

\textbf{(4) MDPs of more genernal single-server queues}

Stidham \cite{Sti:1978}, Crabill \cite{Cra:1972}, Lippman \cite{Lip:1975},
Schassberger \cite{Sch:1975}, Stidham \cite{Sti:1985}, Hordijk and Spieksma
\cite{Hor:1989}, Federgruen and So \cite{Fed:1991}, Lamond \cite{Lam:1991},
Towsley et al. \cite{Tow:1992}, Koole \cite{Koo:1998}, Haviv and Puterman
\cite{Hav:1998}, Lewis et al. \cite{Lew:1999}, George and Harrison
\cite{Geo:2001}, Johansen and Larsen \cite{Joh:2001}, Piunovskiy
\cite{Piu:2004}, Stidham \cite{Sti:2005}, Adusumilli and Hasenbein
\cite{Adu:2010}, Kumar et al. \cite{Kum:2013} and Yan\ et al. \cite{Yan:2013}.

\textbf{(5) MDPs of single-server batch queues}

Deb and Serfozo \cite{Deb:1973}, Deb \cite{Deb:1976} and Powell and Humblet
\cite{Pow:1986} with batch services; and Nobel and Tijms \cite{Nob:1999} with
batch arrivals.

\textbf{(6) MDPs of single-server queues with either balking, reneging or abandonments}

Blackburn \cite{Bla:1972} with balking, Down et al. \cite{Dow:2011} with
reneging, and Legros \cite{Leg:2018} with abandonments.

\textbf{(7) MDPs of single-server priority queues}

Robinson \cite{Rob:1978}, Browne and Yechiali \cite{Bro:1989}, Groenevelt et
al. \cite{Gro:2002} and Brouns and Van Der Wal \cite{Bro:2006}.

\textbf{(8) MDPs of single-server processor-sharing queues}

De Waal \cite{Waa:1993}, Altman et al. \cite{Alt:2001}, Van der Weij et al.
\cite{Van:2008} and Bhulai et al. \cite{Bhu:2014}.

\textbf{(9) MDPs of single-server retrial queues}

Liang and Kulkarni \cite{Lia:1999}, Winkler \cite{Win:2013} and Giovanidis et
al. \cite{Gio:2009}.

\textbf{(10) MDPs of single-server information-based queues}

Kuri and Kumar \cite{Kur:1995, Kur:1997}, Altman and Stidham \cite{Alt:1995}
and Honhon and Seshadri \cite{Hon:2007}.

\textbf{(11) MDPs of single-server queues with multiple classes of customers}

Harrison \cite{Har:1975}, Chen \cite{Che:1989}, Browne and Yechiali
\cite{Bro:1991}, De Serres \cite{De:1991a, De:1991b}, Ata \cite{Ata:2006},
Feinberg and Yang \cite{Fei:2011} and Larra\~{n}aga et al.\ \cite{Lar:2014}.

\textbf{(12) MDPs of single-server queues with optimal pricing}

Low \cite{Low:1974}, Chen \cite{Che:1989}, Yoon and Lewis \cite{Yoo:2004},
\c{C}elik and Maglaras \cite{Cel:2008}, Economou and Kanta \cite{Eco:2008}
and\ Yildirim and Hasenbein \cite{Yil:2010}.

\textbf{(13) MDPs of single-server manufacturing queues}

\textit{(a) The make-to-stock queues: }Sava\c{s}aneril et al. \cite{Sav:2010},
Sanajian et al. \cite{San:2010}, Perez and Zipkin \cite{Per:1997}, Jain
\cite{Jai:2006} and Cao and Xie \cite{Cao:2016}.

\textit{(b) The make-to-order queues: }Besbes and Maglaras \cite{Bes:2009} and
\c{C}elik and Maglaras \cite{Cel:2008}.

\textit{(c) The assemble-type queues:} Nadar et al. \cite{Nad:2014}.

\textit{(d) The inventory control queues:} Veatch \cite{Vea:1992},
Sava\c{s}aneril et al. \cite{Sav:2010}, Federgruen and Zipkin \cite{Fed:1984},
Federgruen and Zheng \cite{Fed:1992}, Feinberg \cite{Fed:2016}, Feinberg and
Liang \cite{Fed:2017}.

\textbf{(14) MDPs of inventory rationing across multiple demand classes}

Ha \cite{HaA:1997, HaB:1997, Ha:2000}, Gayon et al. \cite{Gay:2009} and Li et
al. \cite{Li:2019}.

\section{MDPs of \textbf{Multi-server Queues}}

In this section, we provide an overview for MDPs of multi-server queues, which
are another important research direction.

\textbf{(1) MDPs of M/M/c queues}

Low \cite{Low:1974b}, Anderson \cite{And:1983}, Printezis and Burnetas
\cite{Pri:2008} and Feinberg and Yang \cite{Fei:2014, Fei:2011}.

\textbf{(2) MDPs of GI/M/c queues}

Yechiali \cite{Yec:1972}, Van Nunen and Puterman \cite{Van:1983} and Feinberg
and Yang \cite{Fei:2016a}.

\textbf{(3) MDPs of two-server queues}

Larsen and Agrawala \cite{Lar:1983}, Lin and Kumar \cite{Lin:1984}, Hajek
\cite{Haj:1984}, Varma \cite{Var:1991}, Chen et al. \cite{Che:1994} and Xu and
Zhao \cite{Xu:1996}.

\textbf{(4) MDPs of multi-server queues}

Emmons \cite{Emm:1972}, Helm and Waldmann \cite{Hel:1984}, Blanc et al.
\cite{Bla:1992}, Bradford \cite{Bra:1996}, Ko\c{c}a\v{g}a and Ward
\cite{Koc:2010} and Lee and Kulkarni \cite{Lee:2014}.

\textbf{(5) MDPs of heterogeneous server queues}

Rosberg and Kermani \cite{Ros:1989}, Nobel and Tijms \cite{Nob:2000}, Rykov
\cite{Rky:2001}, Rykov and Efrosinin \cite{Rky:2004} and Tirdad et al.
\cite{Tir:2016}.

\section{MDPs of \textbf{Queueing Networks}}

In this section, we provide an overview for MDPs of queueing networks. Note
that the MDPs of queueing networks have been an interesting research direction
for many years, and they have also established key applications in many
practical areas.

\textbf{(1) MDPs of more general queueing networks}

Ross \cite{Ros:1986}, Weber and Stidham \cite{Web:1987}, Stidham
\cite{Sti:1988}, Shanthikumar and Yao \cite{Sha:1989}, Veatch and Wein
\cite{Vea:1992b}, Tassiulas and Ephremides\cite{Tas:1996}, Papadimitriou and
Tsitsiklis \cite{Pap:1999}, B\"{a}uerle \cite{Bau:2000}, B\"{a}uerle
\cite{Bau:2002} and Solodyannikov \cite{Sol:2014}.

\textbf{(2) MDPs of queueing networks with multiple classes of customers}

Shioyama \cite{Shi:1991}, Bertsimas et al. \cite{Ber:1994}, Maglaras
\cite{Mag:1999}, Chen and Meyn \cite{Che:1999} and Cao and Xie
\cite{Cao:2016b}.

\textbf{(3) Queueing applications of Markov decision processes}

Serfozo \cite{Ser:1981} studied the MDPs of birth-death processes and random
walks, and then discussed dynamic control queueing networks. White
\cite{Whi:2005} focused on the MDPs of QBD processes, which were used to deal
with dynamic control of queueing networks. Robinson \cite{Rob:1976} and
Hordijk et al. \cite{Hor:1983} studied the MDP which were applied to the study
of queueing networks. Sennott \cite{Sen:1989} analyzed the semi-MDP and
applied the obtained results to discuss the queueing networks.

Other key research includes Van Dijk and Puterman \cite{Van:1988}, Liu et al.
\cite{Liu:1995}, Altman et al. \cite{Alt:2003} and Adlakha et al.
\cite{Adl:2012}.

\section{MDPs of \textbf{Queueing Networks with Special Structure}}

In this section, we provide an overview for MDPs of queueing networks with
special Structure, for example, multi-station tandem queues, multi-station
parallel queues, polling queues, fork-join queues and so on.

\textbf{(1) MDPs of two-station tandem queues}

Ghoneim and Stidham \cite{Gho:1985}, Nishimura \cite{Nis:1986}, Farrar
\cite{Far:1993}, Iravani et al. \cite{Ira:1997}, Ahn et al. \cite{Ahn:2002}
and Zayas-Cab\'{a}n et al. \cite{Zay:2016}.

\textbf{(2) MDPs of multi-station tandem queues}

Rosberg et al. \cite{Ros:1982}, Hordijk and Koole \cite{Hor:1992}, Hariharan
et al. \cite{Har:1996}, Gajrat et al. \cite{Gaj:2003}, Koole \cite{Koo:2004},
Zhang and Ayhan \cite{Zha:2013} and Leeuwen and N\'{u}nez-Queija
\cite{Lee:2017}.

\textbf{(3) MDPs of parallel queues}

parallel queues by Weber \cite{Web:1978}, Bonomi \cite{Bon:1990}, Menich and
Serfozo \cite{Men:1991}, Xu et al. \cite{Xu:1992}, Hordijk and Koole
\cite{Hor:1993}, Chen et al. \cite{Cha:1994}, Sparaggis et al. \cite{Spa:1996}%
, Koole \cite{Koo:1996}, Ku and Jordan \cite{Ku:2002}, Down and Lewis
\cite{Dow:2006}, Delasay et al. \cite{Del:2012} and Feinberg and Zhang
\cite{Fei:2015}.

\textbf{(4) MDPs of polling queues}

Browne and Yechiali \cite{Bro:1989b}, Gandhi and Cassandras \cite{Gan:1996},
Koole and Nain \cite{Koo:2000} and Gaujal et al. \cite{Gau:2007}.

\textbf{(5) MDPs of fork-Join queueing networks}

Pascual et al. \cite{Pas:2013}, Zeng et al. \cite{Zen:2016}, Marin and Rossi
\cite{Mar:2017} and Zeng et al. \cite{Zen:2018}.

\textbf{(6) MDPs of Call Centers}

Koole \cite{Koo:2013}, Bhulai \cite{Bhu:2009}, Legros et al. \cite{Leg:2016},
Gans et al. \cite{Gan:2003} and Koole and Mandelbaum \cite{Koo:2002}.

\textbf{(7) MDPs of distributed queueing networks}

Chou and Abraham \cite{Cho:1982}, e Silva and Gerla \cite{Sil:1991}, Franken
and Haverkort \cite{Fra:1996}, Li and Kameda \cite{Li:1998}, Nadar et al.
\cite{Nad:2014} and Vercraene et al. \cite{Ver:2018}.

\textbf{(8) Competitive MDPs of distributed queueing networks}

The competitive MDPs are called to be stochastic games. Altman and Hordijk
\cite{Alt:1995b} studied the zero-sum Markov game and applied the obtained
results to the worst-case optimal control of queueing networks. Altman
\cite{Alt:1996b} studied non-zero stochastic games and applied their results
to admission, service and routing control in queueing networks. Altman
\cite{Alt:1999b} proposed a Markov game approach for analyzing the optimal
routing of a queueing network. Hordijk et al. \cite{Hor:1997} studied a
multi-chain stochastic game which was applied to the worst case admission
control in a queueing network. Xu and Hajek \cite{Xu:2013} studied the game
problem of supermarket models. Xia \cite{Xia:2014} applied the stochastic
games to analyzing the service rate control of a closed queueing network.

\textbf{(9) Heavy traffic analysis for controlled queues and networks }

Heavy traffic analysis can be used to deal with a class of important problems
of controlled queues and networks by means of fluid and diffusion
approximation. Readers may refer to, for example, Kushner \cite{Kus:2001},
Kushner and Ramachandran \cite{Kus:1989}, Kushner and Martins \cite{Kus:1996};
Harrison \cite{Har:1985}, Plambeck et al. \cite{Pla:2001}; Chen and Yao
\cite{ChenY:2001}, Atar et al. \cite{Ata:2004}.

\section{Key Objectives in Queueing Dynamic Control}

In this section, we introduce some key objectives to classify the literature
of queueing dynamic control, for example, input control, service control,
dynamic control under different service mechanisms, dynamic control with
pricing, threshold control and so forth.

\textbf{Objective one: Input control}

The input control is to apply the MDPs to dynamically control the input
process of customers in the queues and networks, including the input rate
control, the interval time control, and the admission access control (e.g.,
probability that an arriving customer chooses entering the system or some servers).

\textit{(a) The input rate control: }Kitaev and Rykov \cite{Kit:1995}, Sennott
\cite{Sen:2009}, Crabill et al. \cite{Cra:1977}, Stidham and Weber
\cite{Sti:1993}, Crabill \cite{Cra:1972} and Lee and Kulkarni \cite{Lee:2014}.

\textit{(b) The input process control:} Kitaev and Rykov \cite{Kit:1995},
Sennott \cite{Sen:2009}, Crabill et al. \cite{Cra:1977}, Stidham and Weber
\cite{Sti:1993}, Abdel-Gawad \cite{Abd:1984}, Doshi \cite{Dos:1977}, Stidham
\cite{Sti:1978}, Piunovskiy \cite{Piu:2004}, Kuri and Kumar \cite{Kur:1995},
Kuri and Kumar \cite{Kur:1997}, Van Nunen and Puterman \cite{Van:1983}, Helm
and Waldmann \cite{Hel:1984}, Ghoneim and Stidham \cite{Gho:1985} and
Nishimura \cite{Nis:1986}.

\textit{(c) The admission access control: }Crabill et al. \cite{Cra:1973,
Cra:1977}, Stidham and Weber \cite{Sti:1993}, Brouns \cite{Bro:2003}, Rue and
Rosenshine \cite{Rue:1981a, Rue:1981b, Rue:1983}, Dimitrakopoulos and Burnetas
\cite{Dim:2017}, Mandelbaum and Yechiali \cite{Man:1983}, Mendelson and
Yechiali \cite{Men:1981}, Stidham \cite{Sti:1985}, Hordijk and Spieksma
\cite{Hor:1989}, Lamond \cite{Lam:1991}, Lewis et al. \cite{Lew:1999},
Adusumilli and Hasenbein \cite{Adu:2010}, Altman et al. \cite{Alt:2001},
Honhon and Seshadri \cite{Hon:2007}, Yoon and Lewis \cite{Yoo:2004}, Yildirim
and Hasenbein \cite{Yil:2010}, Anderson \cite{And:1983}, Emmons
\cite{Emm:1972}, Blanc et al. \cite{Bla:1992}, Ko\c{c}a\v{g}a and Ward
\cite{Koc:2010}, Zhang and Ayhan \cite{Zha:2013}, Altman \cite{Alt:1996b},
Hordijk et al. \cite{Hor:1997} and Xia \cite{Xia:2014}.

\textbf{Objective two: Service control}

The service control is to use the MDPs to dynamically control the service
process in queues and networks, including the service rate control, the
service time control, and the service process control.

\textit{(a) The service rate control: }Kitaev and Rykov \cite{Kit:1995},
Sennott \cite{Sen:2009}, Stidham \cite{Sti:2002, Sti:2009}, Crabill et al.
\cite{Cra:1973, Cra:1977}, Stidham and Weber \cite{Sti:1993}, Yao and
Schechner \cite{Yao:1989}, Dimitrakopoulos and Burnetas \cite{Dim:2017},
Mitchell \cite{Mit:1973}, Doshi \cite{Dos:1978}, Jo and Stidham \cite{Jo:1983}%
, Adusumilli and Hasenbein \cite{Adu:2010}, Kumar et al. \cite{Kum:2013},
Anderson \cite{And:1983}, Lee and Kulkarni \cite{Lee:2014}, Weber and Stidham
\cite{Web:1987}, Ma and Cao \cite{Ma:1994}, Xia \cite{Xia:2014}, Xia and
Shihada \cite{Xia:2013} and Xia and Jia \cite{Xia:2015}.

\textit{(b) The service time control: }Gallisch \cite{Gal:1979}.

\textit{(c) The service process control: }Kitaev and Rykov \cite{Kit:1995},
Sennott \cite{Sen:2009}, Stidham \cite{Sti:2009}, Crabill et al.
\cite{Cra:1973, Cra:1977}, Stidham and Weber \cite{Sti:1993}, Schassberger
\cite{Sch:1975}, Johansen and Larsen \cite{Joh:2001}, Stidham \cite{Sti:2005},
Nishimura \cite{Nis:1986}, Rosberg et al. \cite{Ros:1982}, Altman
\cite{Alt:1996b} and Hordijk et al. \cite{Hor:1997}.

\textbf{Objective three: Dynamic control under different service mechanisms}

Many practical and real problems lead to introduction of different service
mechanisms which make some interesting queueing systems, for example, priority
queues, processor-sharing queues, retrial queues, vacation queues, repairable
queues, fluid queues and so on.

\textit{(a) The priority queueing control: }The priority is an important
service mechanism, and it is a precondition that sets up useful relations
among key customers, segmenting market and adhering to long-term cooperation.
Note that the priority makes dynamic control of queues with multi-class
customers. Readers may refer to Rykov and Lembert \cite{Ryk:1967}, Crabill et
al. \cite{Cra:1973, Cra:1977}, Stidham and Weber \cite{Sti:1993}, Kofman and
Lippman \cite{Kof:1981}, Robinson \cite{Rob:1978}, Browne and
Yechiali\cite{Bro:1989}, Groenevelt et al. \cite{Gro:2002}, Brouns and Van Der
Wal \cite{Bro:2006}, Jain \cite{Jai:2006}, Printezis and
Burnetas\cite{Pri:2008} and Koole and Nain \cite{Koo:2000}.

\textit{(b) The processor-sharing queueing control: }Crabill et al.
\cite{Cra:1973, Cra:1977}, Stidham and Weber \cite{Sti:1993}, De Waal
\cite{Waa:1993}, Altman et al. \cite{Alt:2001}, Van der Weij et al.
\cite{Van:2008}, Bhulai et al. \cite{Bhu:2014} and Bonomi \cite{Bon:1990}.

\textit{(c) The retrial queueing control: }Bhulai et al. \cite{Bhu:2014},
Liang and Kulkarni \cite{Lia:1999}, Winkler \cite{Win:2013} and Giovanidis et
al. \cite{Gio:2009}.

\textit{(d) The vacation queueing control: }Li et al. \cite{Li:2014}, Altman
and Nain \cite{Alt:1993, Alt:1996}, Kella \cite{Kel:1990} and Federgruen and
So \cite{Fed:1991}.

\textit{(e) The repairable queueing control: }Dimitrakos and Kyriakidis
\cite{Dim:2008}, Rykov and Efrosinin \cite{Ryk:2012}, Tijms and van der Duyn
Schouten \cite{Tij:1985}.

\textit{(f) The removable server control: }For dynamic control of working
servers, it is necessary to real-time response to the peak period or an
emergency phenomenon through increasing or decreasing the number of working
servers according to either customer number or system workload. We refer the
readers to Feinberg and Kim \cite{Fei:1996}, Feinberg and Kella
\cite{Fei:2002} and Iravani et al. \cite{Ira:1997}.

\textit{(g) The dynamic control of queueing behavior:} blocking by Blackburn
\cite{Bla:1972} and Economou and Kanta \cite{Eco:2008}; reneging and
impatience by Li et al. \cite{Li:2014} and Anderson \cite{And:1983}; and
abandonment by Down et al. \cite{Dow:2011}, Legros et al. \cite{Leg:2018},
Larra\~{n}aga et al. \cite{Lar:2014} and Zayas-Cab\'{a}n et al.
\cite{Zay:2016}.

\textbf{Objective four: Threshold control}

In dynamic control of queues and networks, the threshold-type policy is a
simple and effective mode, including single-threshold and dual-threshold.

\textit{(a) The single-threshold policy:} Brouns \cite{Bro:2003}, Altman and
Nain \cite{Alt:1996}, Federgruen and So \cite{Fed:1991}, Brouns and Van Der
Wal \cite{Bro:2006};

\textit{(b) The dual-threshold policy:} Lu and Serfozo \cite{Lu:1984}, Plum
\cite{Plu:1991} and Kitaev and Serfozo \cite{Kit:1999}.

\textbf{Objective five: Optimal routing control}

\textit{(a) The entering parallel-server policy:} Rosberg et al.
\cite{Ros:1982}, Weber \cite{Web:1978}, Bonomi \cite{Bon:1990}, Menich and
Serfozo \cite{Men:1991}, Xu et al. \cite{Xu:1992}, Hordijk and Koole
\cite{Hor:1993}, Chang et al. \cite{Cha:1994}, Sparaggis et al.
\cite{Spa:1996}, Koole \cite{Koo:1996}, Ku and Jordan \cite{Ku:2002}, Down and
Lewis \cite{Dow:2006}, Delasay et al. \cite{Del:2012} and Li and Kameda
\cite{Li:1998}.

\textit{(b) The routing policy:} Abdel-Gawad \cite{Abd:1984}, Altman
\cite{Alt:1994, Alt:1999}, Towsley et al. \cite{Tow:1992}, Liang and Kulkarni
\cite{Lia:1999}, Xu and Zhao \cite{Xu:1996}, Bradford \cite{Bra:1996}, Rosberg
and Kermani \cite{Ros:1989}, Ross \cite{Ros:1986}, Stidham \cite{Sti:1988},
Tassiulas and Ephremides \cite{Tas:1996}, Menich and Serfozo \cite{Men:1991},
Koole \cite{Koo:1996}, Browne and Yechiali \cite{Bro:1989}, Altman and Nain
\cite{Alt:1996} and Ho and Cao \cite{Ho:1985}.

\textit{(c) The assignment policy:} Weber \cite{Web:1978}, Bonomi
\cite{Bon:1990} and Xu et al. \cite{Xu:1992}.

\textbf{Objective six: Controlled queues and networks with useful information}

In the queueing networks, the useful information plays a key role in dynamic
control of queueing networks. Readers may refer to Kuri and Kumar
\cite{Kur:1995}, Altman and Stidham \cite{Alt:1995}, Honhon and Seshadri
\cite{Hon:2007}, Altman et al. \cite{Alt:2004}, Altman and Jim\'{e}nez
\cite{Alt:2013} and Rosberg and Kermani \cite{Ros:1989}.

Load balancing is an interesting research direction in queueing networks with
simply observable information, e.g., see Down and Lewis \cite{Dow:2006}, Chou
and Abraham \cite{Cho:1982}, e Silva and Gerla \cite{Sil:1991}, Li and Kameda
\cite{Li:1998}, Li et al. \cite{Li:2014, Li:2015}, Li \cite{Li:2016} and Li and Lui
\cite{Lilui:2016}.

\textbf{Objective seven: Controlled queues and networks with pricing}

The optimal pricing policy is an important research direction in dynamic
control of queues and networks, e.g., see literature Low \cite{Low:1974}, Chen
and Frank \cite{Che:2001}, Yoon and Lewis \cite{Yoo:2004}, \c{C}elik and
Maglaras \cite{Cel:2008}, Economou and Kanta \cite{Eco:2008}, Yildirim and
Hasenbein \cite{Yil:2010}, Feinberg and Yang \cite{Fei:2016a}, Bradford
\cite{Bra:1996}, Xia and Chen \cite{Xia:2018} and Federgruen and Zheng
\cite{Fed:1992}.

\section{Sensitivity-Based Optimization for MDPs of Queueing Networks}

In this section, we simple introduce the sensitivity-based optimization in the
MDPs, and then provide an overview on how to apply the sensitivity-based
optimization in dynamic control of queues and networks.

In the late 1980s, to study dynamic control of queueing systems, Professors
Yu-Chi Ho and Xi-Ren Cao proposed and developed the infinitesimal perturbation
method for discrete event dynamic systems (DEDS), which is a new research
direction for online simulation optimization of the DEDS. See \cite{Ho:1991}
for more interpretation. Further excellent books include
Glasserman\cite{Gla:1991}, Cao \cite{Cao:1994} and Cassandras and Lafortune
\cite{Cas:2008}.

\textbf{Sensitivity-Based Optimization: }Cao et al. \cite{Cao:1996} and Cao
and Chen \cite{Cao:1997} published a pioneer work that transforms the
infinitesimal perturbation of DEDS, together with the MDPs, into the so-called
sensitivity-based optimization by means of the policy-based Markov processes
and the associated Poisson equations, in which they also developed new
concepts, for example, performance potential, and performance difference
equation. On this research line, Cao \cite{Cao:2007} summarized many basic
results of the sensitivity-based optimization. In addition, Li and Liu
\cite{Li:2004} and Chapter 11 in Li \cite{Li:2010} extended and generalized
the sensitivity-based optimization to a more general perturbed Markov process
with infinite states by means of the RG-factorizations.

So far some work has applied the sensitivity-based optimization to deal with
MDPs of queues and networks, e.g., see Xia and Cao \cite{Xia:2012}, Xia and
Shihada \cite{Xia:2013}, Xia \cite{Xia:2014}, Xia and Jia \cite{Xia:2015}, Xia
et al. \cite{Xia:2017} and Xia and Chen \cite{Xia:2018}; Ma et al.
\cite{MaX:2018, MaL:2018} for data centers; and Li et al. \cite{Li:2019} for
inventory rationing control. It is worthwhile to note that the
sensitivity-based optimization of queues and networks can be effectively
supported and developed by means of the matrix-analytic method by Neuts
\cite{Neu:1981, Neu:1989} and the RG-factorizations by Li \cite{Li:2010}. Also
see Ma et al. \cite{MaX:2018, MaL:2018} and Li et al. \cite{Li:2019} for more details.

Recently, Xi-Ren Cao further extended and generalized the sensitivity-based
optimization to the more general case of diffusion processes, called
\textit{relative optimization of continuous-time and continuous-state
stochastic systems} (see Cao \cite{Cao:2019} with a complete draft). Important
examples include Cao \cite{Cao:2015, Cao:2017, Cao:2018, CaoA:2019, CaoB:2019}
and references therein.

\textbf{Event-Based Optimization Approach: }In many practical systems, an
event usually has a specific physical meaning and can mathematically
correspond to a set of state transitions with the same characteristics. In
general, the number of events from change of system states is much smaller
than the state number of the system. Therefore, such an event can be used to
describe an approximate MDP, hence this sets up a new optimal framework,
called event-based optimization. The event-based optimization can directly
capture the future information and the structure nature of the system, which
are reflected in the event to aggregate performance potential. Note only can
the event-based optimization greatly save the calculation, but also it
alleviates the dimensional disaster of a network decision process.

For the event-based optimization, readers may refer to, for example, dynamic
control of queueing systems by Koole \cite{Koo:1998} and Koole \cite{Koo:2000}%
; dynamic control of Markov systems by Cao \cite{Cao:2005}, Cao
\cite{CaoF:2008}, Xia \cite{Xia:2014c} and Jia \cite{Jia:2011}; partially
observable Markov decision processes by Wang and Cao \cite{Wan:2011}; and
admission control of open queueing networks by Xia \cite{Xia:2014b}.

\section{Concluding Remarks}

In this survey, we provide an overview for the MDPs of queues and networks,
including single-server queues, multi-server queues and queueing networks. At
the same time, the overview is also related to some specific objectives, for
example, input control, service control, dynamic control based on different
service mechanisms, dynamic control based on pricing, threshold control and so on.

Along such a line, there are still a number of interesting directions for
potential future research, for example:

$\bullet$ Developing effective and efficient algorithms to find the optimal
polices and to compute the optimal performance measures, and also probably
linking AI and learning algorithms;

$\bullet$ discussing structure properties of the optimal policy in the MDPs of
queueing networks under intelligent environment (for example, IoT, big data,
cloud service, blockchain and AI), and specifically, dealing with
multi-dimensional queueing dynamic control;

$\bullet$ analyzing structure properties of the optimal policy in the MDPs
with either QBD processes, Markov processes of GI/M/1 type or Markov processes
of M/G/1 type, which are well related to various practical stochastic models.

$\bullet$ applying the sensitivity-based optimization and the event-based
optimization to deal with dynamic control of practical stochastic networks,
for example, production and inventory control, manufacturing control,
transportation networks, healthcare, sharing economics, cloud service,
blockchain, service management, energy-efficient management and so forth.

\section*{Acknowledgements}

Quan-Lin Li was supported by the National Natural Science Foundation of China
under grant No. 71671158 and by the Natural Science Foundation of Hebei
province under grant No. G2017203277. Li Xia was supported by the National
Natural Science Foundation of China under grant No. 61573206. The authors thank X.R. Cao and E.A. Feinberg for their valuable comments and suggestions to improve the presentation of this paper.


\begin{thebibliography}{9}                                                                                                %

\bibitem {Abd:1984}E.F. Abdel-Gawad (1984). Optimal control of arrivals and
routing in a network of queues. Ph.D. dissertation, North Carolina State University.

\bibitem {Adl:2012}S. Adlakha, S. Lall and A. Goldsmith (2012). Networked
Markov decision processes with delays. IEEE Transactions on Automatic Control,
57(4), 1013--1018.

\bibitem {Adu:2010}K.M. Adusumilli and J.J. Hasenbein (2010). Dynamic
admission and service rate control of a queue. Queueing Systems, 66(2), 131--154.

\bibitem {Ahm:2005}M.H. Ahmed (2005). Call admission control in wireless
networks: a comprehensive survey. IEEE Communications Surveys \& Tutorials,
7(1), 49--68.

\bibitem {Ahn:2002}H.S. Ahn, I. Duenyas and M.E. Lewis (2002). Optimal control
of a two-stage tandem queuing system with flexible servers. Probability in the
Engineering and Informational Sciences, 16(4), 453--469.

\bibitem {Als:2015}M.A. Alsheikh, D.T. Hoang, D. Niyato, H.P. Tan and S. Lin
(2015). Markov decision processes with applications in wireless sensor
networks: A survey. arXiv preprint arXiv:1501.00644, Pages 1--29.

\bibitem {Alt:1994}E. Altman (1994). A Markov game approach for optimal
routing into a queuing network. Ph.D. dissertation, INRIA (Institut National
de Recherche en Informatique et en Automatique).

\bibitem {Alt:1996b}E. Altman (1996). Non zero-sum stochastic games in
admission, service and routing control in queueing systems. Queueing Systems,
23(1-4), 259--279.

\bibitem {Alt:1999}E. Altman (1999). Constrained Markov decision processes.
CRC Press.

\bibitem {Alt:1999b}E. Altman. (1999). A Markov game approach for optimal
routing into a queuing network. In: Stochastic and Differential Games,
Birkh\"{a}user, Pages 359--375.

\bibitem {Alt:2002}E. Altman (2002). Applications of Markov decision processes
in communication networks. In: Handbook of Markov Decision Processes,
Springer, Pages 489--536.

\bibitem {Alt:2003}E. Altman, B. Gaujal and A. Hordijk (2003). Discrete-Event
Control of Stochastic Networks: Multimodularity and Regularity. Springer.

\bibitem {Alt:1995b}E. Altman and A. Hordijk (1995). Zero-sum Markov games and
worst-case optimal control of queueing systems. Queueing Systems, 21(3-4), 415--447.

\bibitem {Alt:2013}E. Altman and T. Jim\'{e}nez (2013). Admission control to
an M/M/1 queue with partial information. In: International Conference on
Analytical and Stochastic Modeling Techniques and Applications, pp. 12--21.

\bibitem {Alt:2001}E. Altman, T. Jim\'{e}nez and G. Koole (2001). On optimal
call admission control in resource-sharing system. IEEE Transactions on
Communications, 49(9), 1659--1668.

\bibitem {Alt:2004}E. Altman, T. Jim\'{e}nez, R. N\'{u}\={n}ez Queija and U.
Yechiali (2004). Optimal routing among $\cdot$/M/1 queues with partial
information. Stochastic Models, 20(2), 149--171.

\bibitem {Alt:1993}E. Altman and P. Nain (1993). Optimal control of the M/G/1
queue with repeated vacations of the server. IEEE Transactions on Automatic
Control, 38(12), 1766--1775.

\bibitem {Alt:1996}E. Altman and P. Nain (1996). Optimality of a threshold
policy in the M/M/1 queue with repeated vacations. Mathematical Methods of
Operations Research, 44(1), 75--96.

\bibitem {Alt:1995}E. Altman and S. Stidham (1995). Optimality of monotonic
policies for two-action Markovian decision processes, with applications to
control of queues with delayed information. Queueing Systems, 21(3-4), 267--291.

\bibitem {And:1983}M.Q. Anderson (1983). Optimal admission pricing and service
rate control of an M$^{\text{X}}$/M/s queue with reneging. Naval Research
Logistics, 30(2), 261--270.

\bibitem {And:1991}W.J. Anderson (1991). Continuous-time Markov Chains: An
Applications-oriented Approach. Springer.

\bibitem {Arg:2012}N.T. Argon and Y.C. Tsai (2012). Dynamic control of a
flexible server in an assembly-type queue with setup costs. Queueing Systems,
70(3), 233--268.

\bibitem {Asm:2003}S. Asmussen (2003). Applied Probability and Queues. Springer.

\bibitem {Ata:2006}B. Ata (2006). Dynamic control of a multiclass queue with
thin arrival streams. Operations Research, 54(5), 876--892.

\bibitem {Ata:1997}S.T. Atan (1997). Solution methods for controlled queueing
networks. Ph.D. dissertation, Iowa State University.

\bibitem {Ata:2004}R. Atar, A. Mandelbaum and M.I. Reiman (2004). Scheduling
a multi class queue with many exponential servers: Asymptotic optimality in
heavy traffic. The Annals of Applied Probability, 14(3), 1084--1134.

\bibitem {Bal:2001}S. Balsamo, V. de Nitto Person\'{e} and R. Onvural (2001).
Analysis of Queueing Networks with Blocking. Springer.

\bibitem {Bar:1989}M. Bartroli (1989). On the structure of optimal control
policies for networks of queues. Ph.D. dissertation, University of North
Carolina at Chapel Hill.

\bibitem {Bas:1975}F. Baskett, K.M. Chandy, R.R. Muntz and F.G. Palacios
(1975). Open, closed, and mixed networks of queues with different classes of
customers. Journal of ACM, 22(2), 248--260.

\bibitem {Bau:2000}N. B\"{a}uerle (2000). Asymptotic optimality of tracking
policies in stochastic networks. The Annals of Applied Probability, 10(4), 1065--1083.

\bibitem {Bau:2002}N. B\"{a}uerle (2002). Optimal control of queueing
networks: An approach via fluid models. Advances in Applied Probability,
34(2), 313--328.

\bibitem {Bel:1965}R. Bellman and R.E. Kalaba (1965). Dynamic Programming and
Modern Control Theory. Academic Press.

\bibitem {Ber:1996}D.P. Bertsekas (1995). Dynamic Programming and Optimal
Control. Belmont, Massachusetts: Athena Scientific.

\bibitem {BerT:1996}D.P. Bertsekas and J.N. Tsitsiklis (1996). Neuro-dynamic
programming. Athena Scientific.

\bibitem {Ber:1994}D. Bertsimas, I.C. Paschalidis and J.N. Tsitsiklis (1994).
Optimization of multiclass queueing networks: Polyhedral and nonlinear
characterizations of achievable performance. The Annals of Applied
Probability, 4(1), 43--75.

\bibitem {Bes:2009}O. Besbes and C. Maglaras (2009). Revenue optimization for
a make-to-order queue in an uncertain market environment. Operations Research,
57(6), 1438--1450.

\bibitem {Bhu:2009}S. Bhulai (2009). Dynamic routing policies for multiskill
call centers. Probability in the Engineering and Informational Sciences,
23(1), 101--119.

\bibitem {Bhu:2014}S. Bhulai, A.C. Brooms and F.M. Spieksma (2014). On
structural properties of the value function for an unbounded jump Markov
process with an application to a processor sharing retrial queue. Queueing
Systems, 76(4), 425--446.

\bibitem {Bla:1972}J.D. Blackburn (1972). Optimal control of a single-server
queue with balking and reneging. Management Science, 19(3), 297--313.

\bibitem {Bla:1992}J.P.C. Blanc, P.R. de Waal, P. Nain and D. Towsley (1992).
Optimal control of admission to a multiserver queue with two arrival streams.
IEEE Transactions on Automatic Control, 37(6), 785--797.

\bibitem {Bol:2006}G. Bolch, S. Greiner, H. De Meer and K.S. Trivedi (2006).
Queueing Networks and Markov Chains: Modeling and Performance Evaluation with
Computer Science Applications. John Wiley \& Sons.

\bibitem {Bon:1990}F. Bonomi (1990). On job assignment for a parallel system
of processor sharing queues. IEEE Transactions on Computers, 39(7), 858--869.

\bibitem {Bou:2011}R.J. Boucherie and N.M. Van Dijk (Eds.) (2011). Queueing
Networks: A Fundamental Approach. Springer.

\bibitem {Bou:2017}R.J. Boucherie and N.M. Van Dijk (Eds.) (2017). Markov
Decision Processes in Practice. Springer.

\bibitem {Bra:1996}R.M. Bradford (1996). Pricing, routing, and incentive
compatibility in multiserver queues. European Journal of Operational Research,
89(2), 226--236.

\bibitem {Bra:2008}M. Bramson (2008). Stability of Queueing Networks. Springer.

\bibitem {Bro:1989}S. Browne and U. Yechiali (1989). Dynamic priority rules
for cyclic-type queues. Advances in Applied Probability, 21(2), 432--450.

\bibitem {Bro:1989b}S. Browne and U. Yechiali (1989). Dynamic routing in
polling systems. Teletraffic Science, ITC-12, 1455--1466.

\bibitem {Bro:1991}S. Browne and U. Yechiali (1991). Dynamic scheduling in
single-server multiclass service systems with unit buffers. Naval Research
Logistics, 38(3), 383--396.

\bibitem {Bro:2003}G.A. Brouns (2003). Queueing models with admission and
termination control: monotonicity and threshold results. Technische
Universiteit Eindhoven, Pages 1--198.

\bibitem {Bro:2006}G.A. Brouns and J. Van Der Wal (2006). Optimal threshold
policies in a two-class preemptive priority queue with admission and
termination control. Queueing Systems, 54(1), 21--33.

\bibitem {Buz:1993}J.A. Buzacott and J.G. Shanthikumar (1993). Stochastic
Models of Manufacturing Systems. Prentice Hall.

\bibitem {CaoF:2008}F. Cao (2008). Event-based optimization for the
continuous-time Markov systems. Doctoral dissertation, Hong Kong University of
Science and Technology, Hong Kong.

\bibitem {Cao:2016}P. Cao and J. Xie (2016). Optimal control of an inventory
system with joint production and pricing decisions. IEEE Transactions on
Automatic Control, 61(12), 4235--4240.

\bibitem {Cao:2016b}P. Cao and J. Xie (2016). Optimal control of a multiclass
queueing system when customers can change types. Queueing Systems, 82(3-4), 285--313.

\bibitem {Cao:1994}X.R. Cao (1994). Realization Probabilities: The Dynamics of
Queuing Systems. Springer-Verlag.

\bibitem {Cao:2005}X.R. Cao (2005). Basic ideas for event-based optimization
of Markov systems. Discrete Event Dynamic Systems, 15(2), 169--197.

\bibitem {Cao:2007}X.R. Cao (2007). Stochastic Learning and Optimization: A
Sensitivity-Based Approach. Springer.

\bibitem {Cao:2015}X.R. Cao (2015). Optimization of average rewards of time
nonhomogeneous Markov chains. IEEE Transactions on Automatic Control, 60(7), 1841--1856.

\bibitem {Cao:2017}X.R. Cao (2017). Optimality conditions for long-run average
rewards with underselectivity and nonsmooth features. IEEE Transactions on
Automatic Control, 62(9), 4318--4332.

\bibitem {Cao:2018}X.R. Cao (2018). Semismooth potentials of stochastic
systems with degenerate diffusions. IEEE Transactions on Automatic Control,
63(10), 3566--3572.

\bibitem {Cao:2019}X.R. Cao (2019). Relative optimization of continuous-time
and continuous-state stochastic systems. The complete draft of a Cao's new book.

\bibitem {CaoA:2019}X.R. Cao (2019). State classification and multi-class
optimization of continuous-time and continuous-state Markov processes. IEEE
Transactions on Automatic Control, Online Publication: March 18, 2019, Pages 1--14.

\bibitem {CaoB:2019}X.R. Cao (2019). Stochastic control of multi-dimensional
systems with relative optimization. IEEE Transactions on Automatic Control,
Online Publication: June 27, 2019, Pages 1--15.

\bibitem {Cao:1997}X.R. Cao and H.F. Chen (1997). Perturbation realization,
potentials, and sensitivity analysis of Markov processes. IEEE Transactions on
Automatic Control, 42(10), 1382--1393.

\bibitem {Cao:1996}X.R. Cao, X.M. Yuan, and L. Qiu (1996). A single sample
path-based performance sensitivity formula for Markov chains. IEEE
Transactions on Automatic Control, Vol. 41, Pages 1814--1817.

\bibitem {Cas:1998}A.R. Cassandra (1998). Exact and Approximate Algorithms for
Partially Observable Markov Decision Processes. Doctoral Dissertation, Brown
University Providence.

\bibitem {Cas:2008}C.G. Cassandras and S. Lafortune (2008). Introduction to
Discrete Event Systems. Springer.

\bibitem {Cel:2008}S. \c{C}elik and C. Maglaras (2008). Dynamic pricing and
lead-time quotation for a multiclass make-to-order queue. Management Science,
54(6), 1132--1146.

\bibitem {Cha:1994}C.S. Chang, R. Nelson and D.D. Yao (1994). Optimal task
scheduling on distributed parallel processors. Performance Evaluation,
20(1-3), 207--221.

\bibitem {Cha:2000}C.S. Chang (2000). Performance Guarantees in Communication
Networks. Springer.

\bibitem {Cha:1999}X. Chao, M. Miyazawa and M. Pinedo (1999). Queueing
Networks: Customers, Signals and Product Form Solutions. Wiley.

\bibitem {Che:1989}H. Chen (1989). Optimal intensity control of a multi-class
queue. Queueing Systems, 5(4), 281-293.

\bibitem {Che:2001}H. Chen and M.Z. Frank (2001). State dependent pricing with
a queue. IIE Transactions, 33(10), 847--860.

\bibitem {Che:1994}H. Chen, P. Yang and D.D. Yao (1994). Control and
scheduling in a two-station queueing network: Optimal policies and heuristics.
Queueing Systems, 18(3-4), 301--332.

\bibitem {ChenY:2001}H. Chen and D.D. Yao (2001). Fundamentals of Queueing
Networks: Performance, Asymptotics, and Optimization. Springer.

\bibitem {Chen:2004}M. Chen (2004). From Markov Chains to Non-equilibrium
Particle Systems. World Scientific.

\bibitem {Che:1999}R.R. Chen and S. Meyn (1999). Value iteration and
optimization of multiclass queueing networks. Queueing Systems, 32(1-3), 65--97.

\bibitem {Cho:1982}T.C.K. Chou and J.A. Abraham (1982). Load balancing in
distributed systems. IEEE Transactions on Software Engineering, (4), 401--412.

\bibitem {Chu:1967}K.L. Chung (1967). Markov Chains. Springer-Verlag.

\bibitem {Coh:1969}J.W. Cohen (1969). The Single Server Queue. North-Holland
Publishing Company, Amsterdan.

\bibitem {Cra:1972}T.B. Crabill (1972). Optimal control of a service facility
with variable exponential service times and constant arrival rate. Management
Science, 18(9), 560--566.

\bibitem {Cra:1973}T.B. Crabill, D. Gross and M.J. Magazine (1973). A survey
of research on optimal design and control of queues. No. Serial T-280,
Washington DC program in Logistics, George Washington University.

\bibitem {Cra:1977}T.B. Crabill, D. Gross and M.J. Magazine (1977). A
classified bibliography of research on optimal design and control of queues.
Operations Research, 25(2), 219--232.

\bibitem {Dad:2001}H. Daduna (2001). Queueing Networks with Discrete Time
Scale: Explicit Expressions for the Steady State Behavior of Discrete Time
Stochastic Networks. Springer.

\bibitem {Dai:1995}J.G. Dai (1995). On positive Harris recurrence of
multiclass queueing networks: a unified approach via fluid limit models. The
Annals of Applied Probability, 5(1), 49--77.

\bibitem {De:1991a}Y. De Serres (1991). Simultaneous optimization of flow
control and scheduling in a single server queue with two job classes.
Operations Research Letters, 10(2), 103--112.

\bibitem {De:1991b}Y. De Serres (1991). Simultaneous optimization of
flow-control and scheduling in a single server queue with two job classes:
Numerical results and approximation. Computers \& Operations Research, 18(4), 361--378.

\bibitem {Deb:1976}R.K. Deb (1976). Optimal control of batch service queues
with switching costs. Advances in Applied Probability, 8(1), 177--194.

\bibitem {Deb:1973}R.K. Deb and R.F. Serfozo (1973). Optimal control of batch
service queues. Advances in Applied Probability, 5(2), 340--361.

\bibitem {Del:2012}M. Delasay, B. Kolfal and A. Ingolfsson (2012). Maximizing
throughput in finite-source parallel queue systems. European Journal of
Operational Research, 217(3), 554--559.

\bibitem {Dem:2011}H. Demirkan, J.C. Spohrer and V. Krishna (Eds.) (2011).
Service Systems Implementation. Springer.

\bibitem {Waa:1993}P. De Waal (1993). A constrained optimization problem for a
processor sharing queue. Naval Research Logistics, 40(5), 719--731.

\bibitem {Dim:2017}Y. Dimitrakopoulos and A. Burnetas (2017). The value of
service rate flexibility in an M/M/1 queue with admission control. IISE
Transactions, 49(6), 603--621.

\bibitem {Dim:2008}T.D. Dimitrakos and E.G. Kyriakidis (2008). A semi-Markov
decision algorithm for the maintenance of a production system with buffer
capacity and continuous repair times. International Journal of Production
Economics, 111(2), 752--762.

\bibitem {Din:2013}H.T. Dinh, C. Lee, D. Niyato and P. Wang (2013). A survey
of mobile cloud computing: architecture, applications, and approaches.
Wireless communications and mobile computing, 13(18), 1587-1611.

\bibitem {Dis:1985}R.L. Disney and D. K\"{o}nig (1985). Queueing networks: A
survey of their random processes. SIAM Review, 27(3), 335-403.

\bibitem {Dob:1990}R.L. Dobrushin, M.Y. Kelbert, A.N. Rybko and Y.M. Suhov
(1990). Qualitative methods of queueing network theory. In: Stochastic
Cellular Systems: Ergodicity, Memory, Morphogenesis (ed. by RL Dobrushin, VM
Kryukov and AL Toom), University Press, Manchester, pp. 183--224.

\bibitem {Doo:1953}J.L. Doob (1953). Stochastic Processes. John Wiley and Sons.

\bibitem {Dos:1977}B.T. Doshi (1977). Continuous time control of the arrival
process in an M/G/1 queue. Stochastic Processes and Their Applications, 5(3), 265--284.

\bibitem {Dos:1978}B.T. Doshi (1978). Optimal control of the service rate in
an M/G/1 queueing system. Advances in Applied Probability, 10(3), 682--701.

\bibitem {Dow:2011}D. Down, G.M. Koole and M. Lewis (2011). Dynamic control of
a single-server system with abandonments. Queueing Systems, 67(1), 63--90.

\bibitem {Dow:2006}D.G. Down and M.E. Lewis (2006). Dynamic load balancing in
parallel queueing systems: Stability and optimal control. European Journal of
Operational Research, 168(2), 509--519.

\bibitem {Dsh:1995}J.H. Dshalalow (1995). Advances in Queueing Theory,
Methods, and Open Problems. CRC Press.

\bibitem {Dsh:1997}J.H. Dshalalow (1997). Frontiers in Queueing: Models and
Applications in Science and Engineering. CRC press.

\bibitem {Eco:2008}A. Economou and S. Kanta (2008). Optimal balking strategies
and pricing for the single server Markovian queue with compartmented waiting
space. Queueing Systems, 59(3-4), 237.

\bibitem {Efr:2004}D. Efrosinin (2004). Controlled queueing systems with
heterogeneous servers. Ph.D. dissertation, Universit\"{a}tsbibliothek
(University of Trier).

\bibitem {Emm:1972}H. Emmons (1972). The optimal admission policy to a
multiserver queue with finite horizon. Journal of Applied Probability, 9(1), 103--116.

\bibitem {Erl:1909}A.K. Erlang (1909). The theory of probabilities and
telephone conversations. Nyt Tidsskrift for Matematik, 20(B), 33--39.

\bibitem {Eth:2005}S.N. Ethier and T.G. Kurtz (2005). Markov Processes:
Characterization and Convergence. John Wiley \& Sons.

\bibitem {Ett:2000}M. Ettl, G.E. Feigin, G.Y. Lin and D.D. Yao (2000). A
supply network model with base-stock control and service requirements.
Operations Research, 48(2), 216--232.

\bibitem {Far:1992}T.M. Farrar (1992). Resource allocation in systems of
queues. Ph.D. dissertation, University of Cambridge.

\bibitem {Far:1993}T.M. Farrar (1993). Optimal use of an extra server in a two
station tandem queueing network. IEEE Transactions on Automatic Control,
38(8), 1296--1299.

\bibitem {Far:1976}W. Farrell (1976). Optimal switching policies in a
non-homogeneous exponential queueing system. Ph.D. dissertation, University of
California at Los Angeles.

\bibitem {Fed:1991}A. Federgruen and K.C. So (1991). Optimality of threshold
policies in single-server queueing systems with server vacations. Advances in
Applied Probability, 23(2), 388--405.

\bibitem {Fed:1992}A. Federgruen and Y.S. Zheng (1992). An efficient algorithm
for computing an optimal $(r,Q)$ policy in continuous review stochastic
inventory systems. Operations research, 40(4), 808--813.

\bibitem {Fed:1984}A. Federgruen and P. Zipkin (1984). An efficient algorithm
for computing optimal $(s,S)$ policies. Operations Research, 32(6), 1268--1285.

\bibitem {Fei:2016}E.A. Feinberg (2016). Optimality conditions for inventory
control. In: Optimization Challenges in Complex, Networked and Risky Systems.
INFORMS TutORials in Operations Research, 14--45.

\bibitem {Fei:2002}E.A. Feinberg and O. Kella (2002). Optimality of D-policies
for an M/G/1 queue with a removable server. Queueing Systems, 42(4), 355--376.

\bibitem {Fei:1996}E.A. Feinberg and D.J. Kim (1996). Bicriterion optimization
of an M/G/1 queue with a removable server. Probability in the Engineering and
Informational Sciences, 10(1), 57--73.

\bibitem {Fei:2017}E.A. Feinberg and Y. Liang (2017). Structure of optimal
policies to periodic-review inventory models with convex costs and backorders
for all values of discount factors. Annals of Operations Research, 1--17.

\bibitem {FeiS:2002}E.A. Feinberg and A. Shwartz (Eds.). (2002). Handbook of
Markov Decision Processes: Methods and Applications. Springer.

\bibitem {Fei:2011}E.A. Feinberg and F. Yang (2011). Optimality of trunk
reservation for an M/M/k/N queue with several customer types and holding
costs. Probability in the Engineering and Informational Sciences, 25(4), 537--560.

\bibitem {Fei:2014}E.A. Feinber and F. Yang (2014). Dynamic price optimization
for an M/M/k/N queue with several customer types. ACM SIGMETRICS Performance
Evaluation Review, 41(3), 25--27.

\bibitem {Fei:2016a}E.A. Feinberg and F. Yang (2016). Optimal pricing for a
GI/M/k/N queue with several customer types and holding costs. Queueing
Systems, 82(1-2), 103--120.

\bibitem {Fei:2015}E.A. Feinberg and X. Zhang (2015). Optimal switching on and
off the entire service capacity of a parallel queue. Probability in the
Engineering and Informational Sciences, 29(4), 483--506.

\bibitem {Fil:2012}J. Filar and K. Vrieze (2012). Competitive Markov Decision
Processes. Springer.

\bibitem {Fra:1996}L.J. Franken and B.R. Haverkort (1996). Reconfiguring
distributed systems using Markov-decision models. In: Proceedings of the
Workshop on Trends in Distributed Systems, Pages 219--228.

\bibitem {Gaj:2003}A. Gajrat, A. Hordijk and A. Ridder (2003).
Large-deviations analysis of the fluid approximation for a controllable tandem
queue. The Annals of Applied Probability, 13(4), 1423--1448.

\bibitem {Gal:1979}E. Gallisch (1979). On monotone optimal policies in a
queueing model of M/G/1 type with controllable service time distribution.
Advances in Applied Probability, 11(4), 870--887.

\bibitem {Gan:1996}A.D. Gandhi and C.G. Cassandras (1996). Optimal control of
polling models for transportation applications. Mathematical and Computer
Modelling, 23(11-12), 1--23.

\bibitem {Gan:2003}N. Gans, G. Koole and A. Mandelbaum (2003). Telephone call
centers: Tutorial, review, and research prospects. Manufacturing \& Service
Operations Management, 5(2), 79--141.

\bibitem {Gar:2006}M. Garavello and B. Piccoli (2006). Traffic Flow on
Networks. Springfield: American Institute of Mathematical Sciences.

\bibitem {Gas:2011}N. Gast and B. Gaujal (2011). A mean eld approach for
optimization in discrete time. Discrete Event Dynamic Systems, 21(1), 63--101.

\bibitem {Gas:2012}N. Gast, B. Gaujal and J.Y. Le Boudec (2012). Mean eld for
Markov decision processes: from discrete to continuous optimization. IEEE
Transactions on Automatic Control, 57(9), 2266--2280.

\bibitem {Gau:2007}B. Gaujal, A. Hordijk and D. Van Der Laan (2007). On the
optimal open-loop control policy for deterministic and exponential polling
systems. Probability in the Engineering and Informational Sciences, 21(2), 157--187.

\bibitem {Gay:2009}J.P. Gayon, F. De Vericourt and F. Karaesmen (2009). Stock
rationing in an M/E$_{\text{r}}$/1 multi-class make-to-stock queue with
backorders. IIE Transactions, 41(12), 1096--1109.

\bibitem {Gel:1998}E. Gelenbe, G. Pujolle, E. Gelenbe and G. Pujolle (1998).
Introduction to Queueing Networks. Wiley.

\bibitem {Geo:2001}J.M. George and J.M. Harrison (2001). Dynamic control of a
queue with adjustable service rate. Operations Research, 49(5), 720--731.

\bibitem {Gho:1985}H.A. Ghoneim and S. Stidham (1985). Control of arrivals to
two queues in series. European Journal of Operational Research, 21 (3), 399--409.

\bibitem {Gio:2009}A. Giovanidis, G. Wunder and J. B\"{u}hler (2009). Optimal
control of a single queue with retransmissions: Delay-dropping tradeoffs. IEEE
Transactions on Wireless Communications, 8(7), 3736--3746.

\bibitem {Gla:1991}P. Glasserman and Y.C. Ho (1991). Gradient Estimation via
Perturbation Analysis. Springer.

\bibitem {Gla:1994}P. Glasserman and D.D. Yao (1994). Monotone Structure in
Discrete-Event Systems. John Wiley \& Sons.

\bibitem {Gro:2002}R. Groenevelt, G. Koole and P. Nain (2002). On the bias
vector of a two-class preemptive priority queue. Mathematical Methods of
Operations Research, 55(1), 107--120.

\bibitem {Guo:2009}X. Guo and O. Hern\'{a}ndez-Lerma (2009). Continuous-Time
Markov Decision Processes. Springer.

\bibitem {HaA:1997}A.Y. Ha (1997). Inventory rationing in a make-to-stock
production system with several demand classes and lost sales. Management
Science, 43(8), 1093--1103.

\bibitem {HaB:1997}A.Y. Ha (1997). Stock-rationing policy for a make-to-stock
production system with two priority classes and backordering. Naval Research
Logistics, 44(5), 457--472.

\bibitem {Ha:2000}A.Y. Ha (2000). Stock rationing in an M/E$_{\text{k}}$/1
make-to-stock queue. Management Science, 46(1), 77--87.

\bibitem {Haj:1984}B. Hajek (1984). Optimal control of two interacting service
stations. IEEE Transactions on Automatic Control, 29(6), 491--499.

\bibitem {Har:1996}R. Hariharan, M.S. Moustafa and S Stidham (1996).
Scheduling in a multi-class series of queues with deterministic service times.
Queueing Systems, 24(1-4), 83--99.

\bibitem {Har:1975}J.M. Harrison (1975). Dynamic scheduling of a multiclass
queue: Discount optimality. Operations Research, 23(2), 270--282.

\bibitem {Har:1985}J.M. Harrison (1985). Brownian Motion and Stochastic Flow
Systems. Wiley.

\bibitem {Hav:1998}M. Haviv and M.L. Puterman (1998). Bias optimality in
controlled queueing systems. Journal of Applied Probability, 35(1), 136--150.

\bibitem {Hel:1984}W.E. Helm and K.H. Waldmann (1984). Optimal control of
arrivals to multiserver queues in a random environment. Journal of Applied
Probability, 21(3), 602--615.

\bibitem {Her:1996}O. Hern\'{a}dez-Lerma and J.B. Lasserre (1996).
Discrete-Time Markov Control Processes: Basic Optimality Criteria. Springer.

\bibitem {Her:1999}O. Hern\'{a}dez-Lerma and J.B. Lasserre (1999). Further
Topics on Discrete-time Markov Control Processes. Springer.

\bibitem {Ho:1985}Y.C. Ho and X.R. Cao (1985). Performance sensitivity to
routing changes in queuing networks and flexible manufacturing systems using
perturbation analysis. IEEE Journal on Robotics and Automation, 1(4), 165--172.

\bibitem {Ho:1991}Y.C. Ho and X.R. Cao (1991). Perturbation Analysis of
Discrete-Event Dynamic Systems. Kluwer Academic Publisher.

\bibitem {Hon:2007}D. Honhon and S. Seshadri (2007). Admission control with
incomplete information to a finite buffer queue. Probability in the
Engineering and Informational Sciences, 21(1), 19--46.

\bibitem {Hor:1992}A. Hordijk and G. Koole (1992). On the shortest queue
policy for the tandem parallel queue. Probability in the Engineering and
Informational Sciences, 6(1), 63--79.

\bibitem {Hor:1993}A. Hordijk and G. Koole (1993). On the optimality of LEPT
and
$\mu$%
c rules for parallel processors and dependent arrival processes. Advances in
Applied Probability, 25(4), 979--996.

\bibitem {Hor:1997}A. Hordijk, O. Passchier and F. Spieksma (1997). Optimal
service control against worst case admission policies: A multichained
stochastic game. Mathematical Methods of Operations Research, 45(2), 281--301.

\bibitem {Hor:1983}A. Hordijk and F.A.V.D.D. Schouten (1983). Average optimal
policies in Markov decision drift processes with applications to a queueing
and a replacement model. Advances in Applied Probability, 15(2), 274--303.

\bibitem {Hor:1989}A. Hordijk and F. Spieksma (1989). Constrained admission
control to a queueing system. Advances in Applied Probability, 21(2), 409--431.

\bibitem {How:1960}R.A. Howard (1960). Dynamic Programming and Markov
Processes. MIT Press, Cambridge, Massachusetts, USA.

\bibitem {Hu:2007}Q. Hu and W. Yue (2007). Markov Decision Processes with
Their Applications. Springer.

\bibitem {Ira:1997}S.M. Iravani, M.J.M. Posner and J.A. Buzacott (1997). A
two-stage tandem queue attended by a moving server with holding and switching
costs. Queueing Systems, 26(3-4), 203--228.

\bibitem {Jac:1957}J.R. Jackson (1957). Networks of waiting lines. Operations
Research, 5(4), 518--521.

\bibitem {Jac:1963}J.R. Jackson (1963). Jobshop-like queueing systems.
Management Science, 10(1), 131--142.

\bibitem {Jai:2006}A. Jain (2006). Priority and dynamic scheduling in a
make-to-stock queue with hyperexponential demand. Naval Research Logistics,
53(5), 363--382.

\bibitem {Jia:2011}Q.S. Jia (2011). On solving event-based optimization with
average reward over infinite stages. IEEE Transactions on Automatic Control,
56(12), 2912--2917.

\bibitem {Jo:1991}K.Y. Jo and O.Z. Maimon (1991). Optimal dynamic load
distribution in a class of flow-type flexible manufacturing systems. European
Journal of Operational Research, 55(1), 71--81.

\bibitem {Jo:1983}K.Y. Jo and S. Stidham (1983). Optimal service-rate control
of M/G/1 queueing systems using phase methods. Advances in Applied
Probability, 15 (3), 616--637.

\bibitem {Joh:2001}S.G. Johansen and C. Larsen (2001). Computation of a
near-optimal service policy for a single-server queue with homogeneous jobs.
European Journal of Operational Research, 134(3), 648--663.

\bibitem {Kar:1968}S. Karlin (1968). A First Course in Stochastic Processes.
Academic press.

\bibitem {Kar:1981}S. Karlin and H.E. Taylor (1981). A Second Course in
Stochastic Processes. Elsevier.

\bibitem {Kel:1990}O. Kella (1990). Optimal control of the vacation scheme in
an M/G/1 queue. Operations Research, 38(4), 724--728.

\bibitem {Kel:1976}F.P. Kelly (1976). Networks of queues. Advances in Applied
Probability, 8(2), 416--432.

\bibitem {Kel:1979}F.P. Kelly (1979). Reversibility and Stochastic Networks.
Cambridge University Press.

\bibitem {Kel:1991}F.P. Kelly (1991). Loss networks. The Annals of Applied
Probability, 1(3), 319--378.

\bibitem {Kem:1976}J.G. Kemeny, J.L. Snell and A.W. Knapp (1976). Denumerable
Markov Chains: with a Chapter of Markov Random Fields by David Griffeath. Springer.

\bibitem {Kit:1995}M.Y. Kitaev and V.V. Rykov (1995). Controlled Queueing
Systems. CRC press.

\bibitem {Kit:1999}M.Y. Kitaev and R.F. Serfozo (1999). M/M/1 queues with
switching costs and hysteretic optimal control. Operations Research, 47(2), 310--312.

\bibitem {Kle:1975}L. Kleinrock (1975). Queueing Systems, Vol. I: Theory.
Wiley Interscience.

\bibitem {Kle:1976}L. Kleinrock (1976). Queueing Systems, Vol. II: Computer
Applications. Wiley Interscience.

\bibitem {Koc:2010}Y.L. Ko\c{c}a\u{g}a and A.R. Ward (2010). Admission control
for a multi-server queue with abandonment. Queueing Systems, 65(3), 275--323.

\bibitem {Kof:1981}E. Kofman and S.A. Lippman (1981). An M/M/1 dynamic
priority queue with optional promotion. Operations Research, 29(1), 174--188.

\bibitem {Kol:2012}A. Kolobov (2012). Planning with Markov decision processes:
An AI perspective. In: Synthesis Lectures on Artificial Intelligence and
Machine Learning, 6(1), Pages 1--210.

\bibitem {Koo:1996}G. Koole (1996). On the pathwise optimal Bernoulli routing
policy for homogeneous parallel servers. Mathematics of Operations Research,
21(2), 469--476.

\bibitem {Koo:1998}G. Koole (1998). The deviation matrix of the M/M/1/$\infty$
and M/M/1/N queue, with applications to controlled queueing models. In:
Proceedings of the 37th IEEE Conference on Decision and Control, Vol. 1, Pages 56--59.

\bibitem {Koo:2002}G. Koole and A. Mandelbaum (2002). Queueing models of call
centers: An introduction. Annals of Operations Research, 113(1-4), 41--59.

\bibitem {Koo:2004}G. Koole (2004). Convexity in tandem queues. Probability in
the Engineering and Informational Sciences, 18(1), 13--31.

\bibitem {Koo:2007}G. Koole (2007). Monotonicity in Markov Reward and Decision
Chains: Theory and Applications. Foundations and Trends in Stochastic Systems,
1(1), 1--76.

\bibitem {Koo:2013}G. Koole (2013). Call Center Optimization. Lulu. com.

\bibitem {Koo:2000}G. Koole and P. Nain (2000). On the value function of a
priority queue with an application to a controlled polling model. Queueing
Systems, 34(1-4), 199--214.

\bibitem {Kri:2016}V. Krishnamurthy (2016). Partially Observed Markov Decision
Processes. Cambridge University Press.

\bibitem {Ku:2002}C.Y. Ku and S. Jordan (2002). Access control of parallel
multiserver loss queues. Performance Evaluation, 50(4), 219--231.

\bibitem {Kum:1990}A. Kumar (1990). Task allocation in multiserver systems---A
survey of results. Sadhana, 15(4-5), 381--395.

\bibitem {Kum:2013}R. Kumar, M.E. Lewis and H. Topaloglu (2013). Dynamic
service rate control for a single-server queue with Markov-modulated arrivals.
Naval Research Logistics, 60(8), 661--677.

\bibitem {Kuo:2006}Y. Kuo (2006). Optimal adaptive control policy for joint
machine maintenance and product quality control. European Journal of
Operational Research, 171(2), 586--597.

\bibitem {Kur:1995}J. Kuri and A. Kumar (1995). Optimal control of arrivals to
queues with delayed queue length information. IEEE Transactions on Automatic
Control, 40(8), 1444--1450.

\bibitem {Kur:1997}J. Kuri and A. Kumar (1997). On the optimal control of
arrivals to a single queue with arbitrary feedback delay. Queueing Systems,
27(1-2), 1--16.

\bibitem {Kus:2001}H. Kushner (2001). Heavy Traffic Analysis of Controlled
Queueing and Communication Networks. Springer.

\bibitem {Kus:1996}H. J. Kushner and L. F. Martins (1996). Heavy traffic
analysis of a controlled multiclass queueing network via weak convergence
methods. SIAM Journal on Control and Optimization, 34(5), 1781--1797.

\bibitem {Kus:1989}H. J. Kushner and K. M. Ramachandran (1989). Optimal and
approximately optimal control policies for queues in heavy traffic. SIAM
Journal on Control and Optimization, 27(6), 1293--1318.

\bibitem {Lak:2013}C. Lakshmi and S.A. Iyer (2013). Application of queueing
theory in health care: A literature review. Operations Research for Health
Care, 2(1-2), 25-39.

\bibitem {Lam:1991}B.F. Lamond (1991). Optimal admission policies for a finite
queue with bursty arrivals. Annals of Operations Research, 28(1), 243-260.

\bibitem {Lar:1983}R.L. Larsen and A.K. Agrawala (1983). Control of a
heterogeneous two-server exponential queueing system. IEEE Transactions on
Software Engineering, (4), 522--526.

\bibitem {Lar:2014}M. Larra\~{n}aga, U. Ayesta and I.M. Verloop (2014). Index
policies for a multi-class queue with convex holding cost and abandonments.
ACM SIGMETRICS Performance Evaluation Review, 42(1), 125--137.

\bibitem {Lat:1999}G. Latouche and V. Ramaswami (1999). Introduction to Matrix
Analytic Methods in Stochastic Modeling. SIAM.

\bibitem {Lau:1999}C.J. Lautenbacher and S. Stidham (1999). The underlying
Markov decision process in the single-leg airline yield-management problem.
Transportation Science, 33(2), 136--146.

\bibitem {Lee:2014}N. Lee and V.G. Kulkarni (2014). Optimal arrival rate and
service rate control of multi-server queues. Queueing Systems, 76(1), 37--50.

\bibitem {Lee:2017}D.V. Leeuwen and R. N\'{u}nez-Queija (2017). Near-optimal
switching strategies for a tandem queue. In: Markov Decision Processes in
Practice, Springer, Pages 439--459.

\bibitem {Leg:2016}B. Legros, O. Jouini and G. Koole (2016). Optimal
scheduling in call centers with a callback option. Performance Evaluation, 95, 1--40.

\bibitem {Leg:2018}B. Legros, O. Jouini and G. Koole (2018). A uniformization
approach for the dynamic control of queueing systems with abandonments.
Operations Research, 66(1), 200--209.

\bibitem {Lew:1999}M.E. Lewis, H. Ayhan and R.D. Foley (1999). Bias optimality
in a queue with admission control. Probability in the Engineering and
Informational Sciences, 13(3), 309--327.

\bibitem {Li:1998}J. Li and H. Kameda (1998). Load balancing problems for
multiclass jobs in distributed/parallel computer systems. IEEE Transactions on
Computers, 47(3), 322--332.

\bibitem {Li:2010}Q.L. Li (2010). Constructive Computation in Stochastic
Models with Applications: The RG-Factorizations. Springer and Tsinghua Press.

\bibitem {Li:2016}Q.L. Li (2016). Nonlinear Markov processes in big networks.
Special Matrices, 4(1), 202--217.

\bibitem {Li:2014}Q.L. Li, G. Dai, J.C.S. Lui and Y. Wang (2014). The
mean-field computation in a supermarket model with server multiple vacations.
Discrete Event Dynamic Systems. 24(4), 473--522.

\bibitem {Li:2015}Q.L. Li, Y. Du, G. Dai and M. Wang (2015). On a doubly
dynamically controlled supermarket model with impatient customers. Computers
\& Operations Research, 55, 76--87.

\bibitem {Li:2019}Q.L. Li, Y.M. Li, J.Y. Ma and H.L. Liu (2019). A complete
algebraic solution for optimal dynamic policy in inventory rationing across
multiple demand classes. Online Publication: arXiv: *******

\bibitem {Li:2004}Q.L. Li and L.M. Liu (2004). An algorithmic approach on
sensitivity analysis of perturbed QBD processes. Queueing Systems, 48(3-4), 365--397.

\bibitem {Lilui:2016}Q.L. Li and J.C.S. Lui (2016). Block-structured supermarket models. Discrete Event Dynamic Systems, 26(2), 147¨C-182.

\bibitem {Lia:1999}H.M. Liang and V.G. Kulkarni (1999). Optimal routing
control in retrial queues. In: Applied Probability and Stochastic Processes,
Springer, Pages 203--218.

\bibitem {Lin:1984}W. Lin and P. Kumar (1984). Optimal control of a queueing
system with two heterogeneous servers. IEEE Transactions on Automatic Control,
29(8), 696--703.

\bibitem {Lip:1975}S.A. Lippman (1975). Applying a new device in the
optimization of exponential queuing systems. Operations Research, 23(4), 687--710.

\bibitem {Liu:1995}Z. Liu, P. Nain and D. Towsley (1995). Sample path methods
in the control of queues. Queueing Systems, 21(3-4), 293--335.

\bibitem {Low:1974}D.W. Low (1974). Optimal pricing for an unbounded queue.
IBM Journal of Research and Development, 18(4), 290--302.

\bibitem {Low:1974b}D.W. Low (1974). Optimal dynamic pricing policies for an
M/M/s queue. Operations Research, 22(3), 545--561.

\bibitem {Lu:1984}F.V. Lu, and R.F. Serfozo (1984). M/M/1 queueing decision
processes with monotone hysteretic optimal policies. Operations Research,
32(5), 1116--1132.

\bibitem {Ma:1994}D.J. Ma and X.R. Cao (1994). A direct approach to
decentralized control of service rates in a closed Jackson network. IEEE
Transactions on Automatic Control, 39(7), 1460--1463.

\bibitem {MaX:2018}J.Y. Ma, L. Xia and Q.L. Li (2018). Optimal
energy-efficient policies for data centers through sensitivity-based
optimization. Online Publication: arXiv: 1808.07905, Pages 1--50.

\bibitem {MaL:2018}J.Y. Ma, Q.L. Li and L. Xia (2019). Optimal asynchronous
dynamic policies in energy-efficient data centers. Online Publication: arXiv:
1901.03371, Pages 1--63.

\bibitem {Mag:1999}C. Maglaras (1999). Dynamic scheduling in multiclass
queueing networks: Stability under discrete-review policies. Queueing Systems,
31(3-4), 171--206.

\bibitem {Man:1983}A. Mandelbaum and U. Yechiali (1983). Optimal entering
rules for a customer with wait option at an M/G/1 queue. Management Science,
29(2), 174--187.

\bibitem {Mar:2017}A. Marin and S. Rossi (2017). Power control in saturated
fork-join queueing systems. Performance Evaluation, 116, 101--118.

\bibitem {Mar:1906}A.A. Markov (1906). Rasprostranenie zakona bol'shih chisel
na velichiny, zavisyaschie drug ot druga. In: Izvestiya
Fiziko-matematicheskogo obschestva pri Kazanskom universitete, 2-ya seriya,
tom 15, Pages 135--156.

\bibitem {Men:1981}H. Mendelson and U. Yechiali (1981). Controlling the GI/M/1
queue by conditional acceptance of customers. European Journal of Operational
Research, 7(1), 77--85.

\bibitem {Men:1991}R. Menich and R.F. Serfozo (1991). Optimality of routing
and servicing in dependent parallel processing systems. Queueing Systems,
9(4), 403--418.

\bibitem {Mey:1996}S.P. Meyn and R.L. Tweedie (1996). Markov Chains and
Stochastic Stability. Springer.

\bibitem {Mil:1967}B.L. Miller (1967). Finite state continuous time Markov
decision processes with applications to a class of optimization problems in
queueing theory. Ph.D. dissertation, Stanford University, California, USA.

\bibitem {Mit:1973}B. Mitchell (1973). Optimal service-rate selection in an
M/G/1 Queue. SIAM Journal on Applied Mathematics, 24(1), 19--35.

\bibitem {Nad:2014}E. Nadar, M. Akan and A. Scheller-Wolf (2014). Technical
note---optimal structural results for assemble-to-order generalized M-systems.
Operations Research, 62(3), 571--579.

\bibitem {Neu:1981}M.F. Neuts (1981). Matrix-Geometric Solutions in Stochastic
Models: An Algorithmic Approach. The Johns Hopkins University Press, Baltimore.

\bibitem {Neu:1989}M.F. Neuts (1989). Structured Stochastic Matrices of M/G/1
Type and Their Applications. Dekker.

\bibitem {Nis:1986}S. Nishimura (1986). Service mechanism control and arrival
control of a two-station tandem queue. Journal of the Operations Research
Society of Japan, 29(3), 191--205.

\bibitem {Nob:1999}R.D. Nobel and H.C. Tijms (1999). Optimal control for an
M$^{\text{X}}$/G/1 queue with two service modes. European Journal of
Operational Research, 113(3), 610--619.

\bibitem {Nob:2000}R.D. Nobel and H.C. Tijms (2000). Optimal control of a
queueing system with heterogeneous servers and setup costs. IEEE Transactions
on Automatic Control, 45(4), 780--784.

\bibitem {Oka:2015}H. Okamura, S. Miyata and T. Dohi (2015). A Markov decision
process approach to dynamic power management in a cluster system. IEEE Access,
3, 3039--3047.

\bibitem {Paj:2014}J. Pajarinen, A. Hottinen and J. Peltonen (2014).
Optimizing spatial and temporal reuse in wireless networks by decentralized
partially observable Markov decision processes. IEEE Transactions on Mobile
Computing, 13(4), 866--879.

\bibitem {Pap:1999}C.H. Papadimitriou and J.N. Tsitsiklis (1999). The
complexity of optimal queuing network control. Mathematics of Operations
Research, 24(2), 293--305.

\bibitem {Pas:2013}R. Pascual, A. Mart\'{\i}nez and R. Giesen (2013). Joint
optimization of fleet size and maintenance capacity in a fork-join cyclical
transportation system. Journal of the Operational Research Society, 64(7), 982--994.

\bibitem {Pat:2011}J. Patrick and M.A. Begen (2011). Markov decision processes
and its applications in healthcare. Handbook of healthcare delivery systems. CRC.

\bibitem {Per:1997}A.P. Perez and P. Zipkin (1997). Dynamic scheduling rules
for a multiproduct make-to-stock queue. Operations Research, 45(6), 919--930.

\bibitem {Piu:2004}A.B. Piunovskiy (2004). Bicriteria optimization of a queue
with a controlled input stream. Queueing Systems, 48(1-2), 159--184.

\bibitem {Pla:2001}E. Plambeck, S. Kumar and J.M. Harrison (2001). A
multiclass queue in heavy traffic with throughput time constraints:
Asymptotically optimal dynamic controls. Queueing Systems, 39(1), 23--54.

\bibitem {Plu:1991}H.J. Plum (1991). Optimal monotone hysteretic Markov
policies in an M/M/1 queueing model with switching costs and finite time
horizon. Zeitschrift f\"{u}r Operations Research, 35(5), 377--399.

\bibitem {Pow:1986}W.B. Powell and P. Humblet (1986). The bulk service queue
with a general control strategy: theoretical analysis and a new computational
procedure. Operations Research, 34(2), 267--275.

\bibitem {Pri:2008}A. Printezis and A. Burnetas (2008). Priority option
pricing in an M/M/m queue. Operations Research Letters, 36(6), 700--704.

\bibitem {Pur:1994}M.L. Puterman (1994). Markov Decision Processes: Discrete
Stochastic Dynamic Programming. John Wiley \& Sons.

\bibitem {Qiu:1999}Q. Qiu and M. Pedram (1999). Dynamic power management based
on continuous-time Markov decision processes. In: Proceedings of the 36th
annual ACM/IEEE Design Automation Conference, Pages 555--561.

\bibitem {Rob:1976}D.R. Robinson (1976). Markov decision chains with unbounded
costs and applications to the control of queues. Advances in Applied
Probability, 8(1), 159--176.

\bibitem {Rob:1978}D.R. Robinson (1978). Optimization of priority queues---a
semi-Markov decision chain approach. Management Science, 24(5), 545--553.

\bibitem {Ros:1989}Z. Rosberg and P. Kermani (1989). Customer routing to
different servers with complete information. Advances in Applied Probability,
21(4), 861--882.

\bibitem {Ros:1982}Z. Rosberg, P. Varaiya and J. Walrand (1982). Optimal
control of service in tandem queues. IEEE Transactions on Automatic Control,
27(3), 600--610.

\bibitem {Ros:1986}K.W. Ross (1986). Optimal dynamic routing in Markov
queueing networks. Automatica, 22(3), 367--370.

\bibitem {Rue:1981a}R.C. Rue and M. Rosenshine (1981). Optimal control for
entry of many classes of customers to an M/M/1 queue. Naval Research
Logistics, 28(3), 489--495.

\bibitem {Rue:1981b}R.C. Rue and M. Rosenshine (1981). Some properties of
optimal control policies for entry to an M/M/1 queue. Naval Research
Logistics, 28(4), 525--532.

\bibitem {Rue:1983}R.C. Rue and M. Rosenshine (1983). Optimal control of entry
to an M/Ek/1 queue serving several classes of customers. Naval Research
Logistics, 30(2), 217--226.

\bibitem {Rue:1985}R.C. Rue and M. Rosenshine, (1985). The application of
semi-Markov decision processes to queueing of aircraft for landing at an
airport. Transportation Science, 19(2), 154--172.

\bibitem {Ryk:1975}V.V. Rykov (1975). Controllable queueing systems (In
Russian). Itogi Nauki i Tekhniki. Seriya'' Teoriya Veroyatnostei.
Matematicheskaya Statistika. Teoreticheskaya Kibernetika'', Vol. 12, 45--152.
(There is English translation in Journal of Soviet Mathematics)

\bibitem {Rky:2001}V.V. Rykov (2001). Monotone control of queueing systems
with heterogeneous servers. Queueing systems, 37(4), 391--403.

\bibitem {Ryk:2017}V.V. Rykov (2017). Controllable queueing systems: From the
very beginning up to nowadays. Reliability: Theory \& Applications, 12(2
(45)), 39--61.

\bibitem {Rky:2004}V.V. Rykov and D. Efrosinin (2004). Optimal control of
queueing systems with heterogeneous servers. Queueing Systems, 46(3-4), 389--407.

\bibitem {Ryk:2012}V.V. Rykov and D. Efrosinin (2012). On optimal control of
systems on their life time. In: Recent Advances in System Reliability, pp. 307--319.

\bibitem {Ryk:1967}V.V. Rykov and E. Lembert (1967). Optimal dynamic
priorities in single-line queueing systems. Engineering Cybernetics, 5(1), 21--30.

\bibitem {San:2010}N. Sanajian, H. Abouee-Mehrizi and B. Balc\i oglu (2010).
Scheduling policies in the M/G/1 make-to-stock queue. Journal of the
Operational Research Society, 61(1), 115--123.

\bibitem {Sav:2010}S. Sava\c{s}aneril, P.M. Griffin and P. Keskinocak (2010).
Dynamic lead-time quotation for an M/M/1 base-stock inventory queue.
Operations research, 58(2), 383--395.

\bibitem {Sch:1975}R. Schassberger (1975). A note on optimal service selection
in a single server queue. Management Science, 21(11), 1326--1331.

\bibitem {Sen:1989}L.I. Sennott (1989). Average cost semi-Markov decision
processes and the control of queueing systems. Probability in the Engineering
and Informational Sciences, 3(2), 247--272.

\bibitem {Sen:2009}L.I. Sennott (2009). Stochastic Dynamic Programming and the
Control of Queueing Systems. John Wiley \& Sons.

\bibitem {Ser:1981}R.F. Serfozo (1981). Optimal control of random walks, birth
and death processes, and queues. Advances in Applied Probability, 13(1), 61--83.

\bibitem {Ser:1999}R.F. Serfozo (1999). Introduction to Stochastic Networks. Springer.

\bibitem {Sha:1989}J.G. Shanthikumar and D.D. Yao (1989). Stochastic
monotonicity in general queueing networks. Journal of Applied Probability,
26(2), 413--417.

\bibitem {Shi:1991}T. Shioyama (1991). Optimal control of a queuing network
system with two types of customers. European Journal of Operational Research,
52(3), 367--372.

\bibitem {Sig:2013}O. Sigaud and O. Buffet (Eds.) (2013). Markov Decision
Processes in Artificial Intelligence. John Wiley \& Sons.

\bibitem {Sil:1991}E.D.S. e Silva and M. Gerla (1991). Queueing network models
for load balancing in distributed systems. Journal of Parallel and Distributed
Computing, 12(1), 24--38.

\bibitem {Sob:1974}M.J. Sobel (1974). Optimal operation of queues. In:
Mathematical Methods in Queueing Theory, Springer, Pages 231--261.

\bibitem {Sol:2014}Y.V. Solodyannikov (2014). Control and observation for
dynamical queueing networks. I. Automation and Remote Control, 75(3), 422--446.

\bibitem {Spa:1996}P.D. Sparaggis, D. Towsley and C.G. Cassandras (1996).
Optimal control of multiclass parallel service systems. Discrete Event Dynamic
Systems, 6(2), 139--158.

\bibitem {Sti:1970}S. Stidham (1970). On the optimality of single-server
queuing systems. Operations Research, 18 (4), 708--732.

\bibitem {Sti:1978}S. Stidham (1978). Socially and individually optimal
control of arrivals to a GI/M/1 queue. Management Science, 24 (15), 1598--1610.

\bibitem {Sti:1985}S. Stidham (1985). Optimal control of admission to a
queueing system. IEEE Transactions on Automatic Control, 30(8), 705--713.

\bibitem {Sti:1988}S. Stidham (1988). Scheduling, routing, and flow control in
stochastic networks. In: Stochastic Differential Systems, Stochastic Control
Theory and Applications, Springer, Pages 529--561.

\bibitem {Sti:2002}S. Stidham (2002). Analysis, design, and control of
queueing systems. Operations Research, 50(1), 197--216.

\bibitem {Sti:2005}S. Stidham (2005). On the optimality of a full-service
policy for a queueing system with discounted costs. Mathematical Methods of
Operations Research, 62(3), 485--497.

\bibitem {Sti:2009}S. Stidham (2009). Optimal Design of Queueing Systems.
Chapman and Hall/CRC.

\bibitem {Sti:1974}S. Stidham and N.U. Prabhu (1974). Optimal control of
queueing systems. In: Mathematical Methods in Queueing Theory, Springer, Pages 263--294.

\bibitem {Sti:1993}S. Stidham and R. Weber (1993). A survey of Markov decision
models for control of networks of queues. Queueing systems, 13(1-3), 291--314.

\bibitem {Sun:2014}L. Sun, H. Dong, F.K. Hussain, O.K. Hussain and E. Chang
(2014). Cloud service selection: State-of-the-art and future research
directions. Journal of Network and Computer Applications, 45, 134--150.

\bibitem {Sys:1997}R. Syski (1997). A personal view of queueing theory. In:
Frontiers in Queueing: Models and Applications in Science and Engineering,
Pages 3--18. CRC Press.

\bibitem {Tas:1996}L. Tassiulas and A. Ephremides (1996). Throughput
properties of a queueing network with distributed dynamic routing and flow
control. Advances in Applied Probability, 28(1), 285--307.

\bibitem {Tij:1994}H.C. Tijms (1994). Stochastic Models: An Algorithmic
Approach. John Wiley \& Sons.

\bibitem {Tij:1985}H.C. Tijms and F. A. van der Duyn Schouten (1985). A Markov
decision algorithm for optimal inspections and revisions in a maintenance
system with partial information. European Journal of Operational Research,
21(2), 245--253.

\bibitem {Tir:2016}A. Tirdad, W.K. Grassmann and J. Tavakoli (2016). Optimal
policies of $M(t)/M/c/c$ queues with two different levels of servers. European
Journal of Operational Research, 249(3), 1124--1130.

\bibitem {Tow:1992}D. Towsley, P.D. Sparaggis and C.G. Cassandras (1992).
Optimal routing and buffer allocation for a class of finite capacity queueing
systems. IEEE Transactions on Automatic Control, 37(9), 1446--1451.

\bibitem {Van:2010}P.T. Vanberkel, R.J. Boucherie, E.W. Hans, J.L. Hurink and
N. Litvak (2010). A survey of health care models that encompass multiple departments.

\bibitem {Van:2008}W. Van der Weij, S. Bhulai and R. Van der Mei (2008).
Optimal scheduling policies for the limited processor sharing queue. Technical
Report WS2008-5, Department of Mathematics, Vrije University.

\bibitem {Dij:1993}N.M. Van Dijk (1993). Queueing Networks and Product Forms:
A Systems Approach. John Wiley \& Son Limited.

\bibitem {Van:1988}N.M. Van Dijk and M.L. Puterman (1988). Perturbation theory
for Markov reward processes with applications to queueing systems. Advances in
Applied Probability, 20(1), 79--98.

\bibitem {Van:1983}J.A.E.E. Van Nunen and M.L. Puterman (1983). Computing
optimal control limits for GI/M/s queuing systems with controlled arrivals.
Management Science, 29(6), 725--734. International Journal of Health
Management and Information, 1(1), 37--69.

\bibitem {Var:1991}S. Varma (1991). Optimal allocation of customers in a two
server queue with resequencing. IEEE Transactions on Automatic Control,
36(11), 1288--1293.

\bibitem {Vea:1992}M.H. Veatch (1992). Queueing control problems for
production/inventory systems. Ph.D. dissertation, Massachusetts Institute of Technology.

\bibitem {Vea:1992b}M.H. Veatch and L.M. Wein (1992). Monotone control of
queueing networks. Queueing Systems, 12(3-4), 391--408.

\bibitem {Ver:2018}S. Vercraene, J.P. Gayon and F. Karaesmen (2018) Effects of
system parameters on the optimal cost and policy in a class of
multidimensional queueing control problems. Operations Research, 66(1), 150--162.

\bibitem {Wak:1991}K. Wakuta (1991). Optimal control of an M/G/1 queue with
imperfectly observed queue length when the input source is finite. Journal of
Applied Probability, 28(1), 210--220.

\bibitem {Wan:2011}D.X. Wang and X.R. Cao (2011). Event-based optimization for
POMDPs and its application in portfolio management. Proceedings of the 18th
IFAC World Congress, 44(1), 3228--3233.

\bibitem {Web:1978}R.R. Weber (1978). On the optimal assignment of customers
to parallel servers. Journal of Applied Probability, 15(2), 406--413.

\bibitem {Web:1987}R.R. Weber and S. Stidham (1987). Optimal control of
service rates in networks of queues. Advances in Applied Probability, 19(1), 202--218.

\bibitem {Whi:2005}L.B. White (2005). A new policy evaluation algorithm for
Markov decision processes with quasi birth-death structure. Stochastic Models,
21(2-3), 785--797.

\bibitem {Win:2013}A. Winkler (2013). Dynamic scheduling of a single-server
two-class queue with constant retrial policy. Annals of Operations Research,
202(1), 197--210.

\bibitem {Wu:2010}C.H. Wu, J.T. Lin and W.C. Chien (2010). Dynamic production
control in a serial line with process queue time constraint. International
Journal of Production Research, 48(13), 3823--3843.

\bibitem {Xia:2014}L. Xia (2014). Service rate control of closed Jackson
networks from game theoretic perspective. European Journal of Operational
Research, 237(2), 546--554.

\bibitem {Xia:2014b}L. Xia (2014). Event-based optimization of admission
control in open queueing networks. Discrete Event Dynamic Systems, 24(2), 133--151.

\bibitem {Xia:2012}L. Xia and X.R. Cao (2012). Performance optimization of
queueing systems with perturbation realization. European Journal of
Operational Research, 218(2), 293--304.

\bibitem {Xia:2018}L. Xia and S. Chen (2018). Dynamic pricing control for open
queueing networks. IEEE Transactions on Automatic Control, 63(10), 3290--3300.

\bibitem {Xia:2017}L. Xia, Q.M. He and A.S. Alfa (2017). Optimal control of
state-dependent service rates in a MAP/M/1 queue. IEEE Transactions on
Automatic Control, 62(10), 4965--4979.

\bibitem {Xia:2015}L. Xia and Q.S. Jia (2015). Parameterized Markov decision
process and its application to service rate control. Automatica, 54, 29--35.

\bibitem {Xia:2014c}L. Xia, Q.S. Jia and X.R. Cao (2014). A tutorial on
event-based optimization---a new optimization framework. Discrete Event
Dynamic Systems, 24(2), 103--132.

\bibitem {Xia:2013}L. Xia and B. Shihada (2013). Max-Min optimality of service
rate control in closed queueing networks. IEEE Transactions on Automatic
Control, 58(4), 1051--1056.

\bibitem {Xu:1992}S.H. Xu, R. Righter and J.G. Shanthikumar (1992). Optimal
dynamic assignment of customers to heterogeneous servers in parallel.
Operations Research, 40(6), 1126--1138.

\bibitem {Xu:1996}S.H. Xu and Y.Q. Zhao (1996). Dynamic routing and jockeying
controls in a two-station queueing system. Advances in Applied Probability,
28(4), 1201--1226.

\bibitem {Xu:2013}J. Xu and B. Hajek (2013). The supermarket game. Stochastic
Systems, 3(2), 405--441.

\bibitem {Yan:2013}R. Yang, S. Bhulai and R. van der Mei (2013). Structural
properties of the optimal resource allocation policy for single-queue systems.
Annals of Operations Research, 202(1), 211--233.

\bibitem {Yao:1989}D.D. Yao and Z. Schechner (1989). Decentralized control of
service rates in a closed Jackson network. IEEE Transactions on Automatic
Control, 34(2), 236--240.

\bibitem {Yec:1972}U. Yechiali (1972). Customers' optimal joining rules for
the GI/M/s queue. Management Science, 18(7), 434--443.

\bibitem {Yeh:1983}L. Yeh and L.C. Thomas (1983). Adaptive control of M/M/1
queues---continuous-time Markov decision process approach. Journal of Applied
Probability, 20(2), 368--379.

\bibitem {Yil:2010}U. Yildirim and J.J. Hasenbein (2010). Admission control
and pricing in a queue with batch arrivals. Operations Research Letters,
38(5), 427--431.

\bibitem {Yoo:2004}S. Yoon and M.E. Lewis (2004). Optimal pricing and
admission control in a queueing system with periodically varying parameters.
Queueing Systems, 47(3), 177--199.

\bibitem {Zay:2016}G. Zayas-Cab\'{a}n, J. Xie, L.V. Green and M.E. Lewis
(2016). Dynamic control of a tandem system with abandonments. Queueing
Systems, 84(3-4), 279--293.

\bibitem {Zen:2016}Y. Zeng, A. Chaintreau, D. Towsley and C.H. Xia (2016). A
necessary and sufficient condition for throughput scalability of fork and join
networks with blocking. ACM SIGMETRICS Performance Evaluation Review, 44(1), 25--36.

\bibitem {Zen:2018}Y. Zeng, A. Chaintreau, D. Towsley and C.H. Xia (2018).
Throughput Scalability Analysis of Fork-Join Queueing Networks. Operations
Research, 66(6),1728--1743.

\bibitem {Zha:2013}B. Zhang and H. Ayhan (2013). Optimal admission control for
tandem queues with loss. IEEE Transactions on Automatic Control, 58(1), 163--167.
\end{thebibliography}
\end{document}